\documentclass[11pt]{article}

\usepackage{amsmath,amssymb,mathtools}
\usepackage{hyperref}

\newtheorem{theorem}{Theorem}
\newtheorem{as}{Theorem}[section]
\newtheorem{prop}[as]{Proposition}
\newtheorem{lemma}[as]{Lemma}

\newtheorem{corollary}[as]{Corollary}

\newtheorem{conj}{Conjecture}
\newtheorem{prob}{Problem}

\newcommand{\qed}{\hspace*{\fill} \rule{7pt}{7pt}}
\newcommand{\Proof}{\noindent{\bf Proof.}\ \ }
\newcommand{\F}{{\cal F}}
\newcommand{\C}{{\cal C}}
\newcommand{\T}{{\cal T}}
\newcommand{\G}{{\cal G}}
\newcommand{\Q}{{\cal Q}}
\renewcommand{\P}{{\cal P}}
\renewcommand{\S}{{\cal S}}
\newcommand{\X}{{\cal X}}

\begin{document}

\title{Vector clique decompositions}

\author{
Raphael Yuster
\thanks{Department of Mathematics, University of Haifa, Haifa
31905, Israel. Email: raphy@math.haifa.ac.il.
This research was supported by the Israel Science Foundation
(grant No. 1082/16).}
}

\date{}

\maketitle

\setcounter{page}{1}

\begin{abstract}

Let $\F_k$ be the set of all graphs on $k$ vertices.
For a graph $G$, a $k$-decomposition is a set of induced subgraphs of $G$, each of which is isomorphic to an element of $\F_k$, such that each pair of vertices of $G$ is
in exactly one element of the set. It is a fundamental result of Wilson that for all $n=|V(G)|$ sufficiently large, $G$ has a $k$-decomposition
if and only if $G$ is $k$-divisible, namely $k-1$ divides $n-1$ and $\binom{k}{2}$ divides $\binom{n}{2}$.

Let ${\bf v} \in {\mathbb R}^{|\F_k|}$ be indexed by $\F_k$.
For a $k$-decomposition $L$ of $G$, let $\nu_{\bf v}(L) = \sum_{F \in \F_k} {\bf v}_F d_{L,F}$ where $d_{L,F}$ is the fraction of elements of $L$
that are isomorphic to $F$.
Let $\nu_{\bf v}(G) = \max_{L} \nu_{\bf v}(L)$ and $\nu_{\bf v}(n)=\min\{\nu_{\bf v}(G):|V(G)|=n\}$
\footnote{If $n$ is not such that graphs on $n$ vertices have a $k$-decomposition, one can synthetically define $\nu_{\bf v}(n) = \nu_{\bf v}(m)$
where $m < n$ is the largest integer such that graphs on $m$ vertices are $k$-decomposable.}.
It is not difficult to prove that the the sequence $\nu_{\bf v}(n)$ has a limit so let
$\nu_{\bf v} = \lim_{n \rightarrow \infty} \nu_{\bf v}(n)$.
Replacing $k$-decompositions with their fractional relaxations, one obtains the (polynomial time computable) fractional analogue
$\nu_{\bf v}^*(G)$ and the corresponding fractional values $\nu^*_{\bf v}(n)$ and $\nu^*_{\bf v}$.
Our first main result is that for each ${\bf v} \in {\mathbb R}^{|\F_k|}$
$$
\nu_{\bf v} = \nu^*_{\bf v}\;.
$$
Furthermore, there is a polynomial time algorithm that produces a decomposition $L$ of a $k$-decomposable graph
such that $\nu_{\bf v}(L) \ge \nu_{\bf v} - o_n(1)$.

A similar result holds when $\F_k$ is the family of all tournaments on $k$ vertices and when $\F_k$ is the family of all edge-colorings of $K_k$.

We use these results to obtain new and improved bounds on several decomposition results. For example, we prove that every $n$-vertex tournament which is $3$-divisible
(namely $n=1,3 \bmod 6$) has a triangle decomposition in which the number of directed triangles is less than $0.0222n^2(1+o(1))$
and that every $5$-decomposable $n$-vertex graph has a $5$-decomposition in which the fraction of cycles of length $5$ is $o_n(1)$.

\vspace*{5pt}

\end{abstract}

\section{Introduction}

The problem of decomposing a large graph $G$ into pairwise edge-disjoint copies of a given graph $F$ has been extensively studied and dates back to a result of Kirkman from 1847 \cite{kirkman-1847},
who proved that $K_n$ has a $K_3$-decomposition whenever $n \equiv 1, 3 \bmod 6$.
These divisibility requirements are necessary as in any decomposition of a graph into triangles, the degree of each vertex must be even and the number of edges must be divisible by $3$.

More generally, for a graph $G$ to have an $F$-decomposition, it must trivially hold that the $\gcd$ of the degree sequence of $G$, denoted by $gcd(G)$ is divisible by
$gcd(F)$ and that the number of edges of $G$, denoted by $e(G)$, is divisible by $e(F)$.
We therefore say that $G$ is {\em $F$-divisible} if these two necessary conditions hold.

The $F$-decomposition problem for $G=K_n$ was completely solved (for large $n$) by Wilson
\cite{wilson-1972-I,wilson-1972-II,wilson-1975,wilson-1975-III}. He proved that whenever $n$ is sufficiently large and
$K_n$ is $F$-divisible (which simply means that $n-1$ is divisible by $gcd(F)$ and $\binom{n}{2}$ is divisible by $e(F)$), then it has an $F$-decomposition.
Recently, this result has been generalized by Keevash \cite{keevash-2014} to the complete uniform hypergraph setting \cite{keevash-2014}. See also Glock et al. \cite{GKLO-2016} for another proof.

Another equivalent way to state Wilson's Theorem is the following. Suppose $\F$ is the set of all spanning subgraphs of $F$.
An $\F$-decomposition of a graph $G$ is a set of subgraphs of $G$, each of which is isomorphic to an element of $\F$, 
such that each pair of vertices of $G$ is in exactly one element of the set (and if this pair is an edge of $G$, then it is also an edge of that element).
So, Wilson's theorem asserts that for $n \ge n_0(F)$, if $K_n$ is $F$-divisible, then any graph $G$ with $n$ vertices has an $\F$-decomposition.
This equivalent statement, leads, however, to a wider set of questions as clearly, if $G$ has an $\F$-decomposition, it has many (in particular, since any vertex permutation
may lead to a distinct $\F$-decomposition). Thus, we can ask about the quality of the various $\F$-decompositions with respect to the distribution
of the members of $\F$ in it.

More formally, for a vector ${\bf v} \in {\mathbb R}^{|\F|}$ indexed by $\F$
and an $\F$-decomposition $L$ of $G$, let $D_{\bf v}(L) = \sum_{H \in \F} {\bf v}_H |L_H|$ where $L_H$ is the subset of $L$ whose elements are isomorphic to $H$.
It will be slightly more convenient to normalize this quantity by defining $d_{L,H}=|L_H|/|L|$ to be the density of $H$ in $L$ and defining
$\nu_{\bf v}(L) = \sum_{H \in \F} {\bf v}_H d_{L,H}$, observing that $\nu_{\bf v}(L) = D_{\bf v}(L)/|L|$ and that $|L|=\binom{n}{2}/e(F)$.

The quality of the decomposition is thus measured by $\nu_{\bf v}(G) = \max_{L} \nu_{\bf v}(L)$ where the maximum is over all $\F$-decompositions of $G$.
We call $\nu_{\bf v}(G)$ the {\em optimal} $\F$-decomposition of $G$ with respect to ${\bf v}$.
So, when $K_n$ has an $F$-decomposition, $\nu_{\bf v}(G)$ is well-defined for every ${\bf v} \in {\mathbb R}^{|\F|}$ and every graph $G$ with $n$ vertices.

Notice that given $F$, we may consider additional structures underlined by $F$ other than the spanning subgraphs of $F$.
For instance, we may define $\F$ to be all possible orientations of $F$ and
then an $\F$-decomposition is defined for {\em tournaments} $G$ and $\nu_{\bf v}(G)$ is defined analogously.
Likewise, we my define $\F$ to be all edge colorings of $F$ with colors from a given set of colors.
In this case, an an $\F$-decomposition is defined for edge colorings $G$ of $K_n$ and $\nu_{\bf v}(G)$ is defined analogously.

In what follows, we will state our results for the case of $F=K_k$, although our results do carry over quite seamlessly to certain more general $F$.
We prefer this approach as it seems to be the most interesting case (in fact, already for some questions arising in the case $k=3$), yet it captures all details of the general proof and since all our applications involve the case where $F=K_k$. So, let $\F_k$ denote the set of all graphs on $k$ vertices, let $\T_k$ denote the set of all tournaments on $k$ vertices, and for a color set $C$, let $\C_k$ denote the set of all edge colorings of $K_k$ with colors from $C$.
Let $G$ be an $n$-vertex graph.
If $K_n$ is $K_k$-divisible (i.e. if $k-1$ divides $n-1$ and $\binom{k}{2}$ divides $\binom{n}{2}$), then $G$ is called {\em $k$-divisible}.
Similarly, $G$ is {\em $k$-decomposable} if $K_n$ has a $K_k$-decomposition. In these terms, Wilson's Theorem asserts that
for all $n$ sufficiently large, $G$ is $k$-decomposable if and only if it is $k$-divisible.

We first observe that computing $\nu_{\bf v}(G)$ is easy for some vectors, while NP-Hard for some others.
Indeed, let ${\bf v} \in {\mathbb R}^{|\F_k|}$ be a constant vector, all entries equal to $c$. In this case,
once we know that $G$ is $k$-decomposable (which we can determine in polynomial time by Wilson's Theorem), we have that $\nu_{\bf v}(G)=c$ as any $k$-decomposition
has this optimal weight. On the other hand, consider the case $k=3$ and the vector which assigns $K_3$ the weight $1$ and assigns the other graphs on three vertices,
the weight zero. Suppose $G$ is $3$-divisible. Now, if $G$ and its complement each have a $K_3$-decomposition,
then we would have $\nu_{\bf v}(G)=e(G)/\binom{n}{2}$. Otherwise, we would have $\nu_{\bf v}(G) < e(G)/\binom{n}{2}$. But determining whether a graph and its complement have
a $K_3$-decomposition is NP-Complete (see \cite{DT-1992,holyer-1981})\footnote{In fact it is proved that deciding if a graph has a $K_3$-decomposition is NP-Complete, but it is straightforward to reduce this problem to the problem of whether a graph and its complement each have a $K_3$-decomposition.}.

Another minor observation is that if ${\bf w} \in {\mathbb R}^{|\F_k|}$ is obtained from ${\bf v}$ by dilation and translation with a constant vector,
namely, ${\bf w} = c{\bf v}+d{\bf 1}$ for some $c > 0$, then $\nu_{\bf w}(G) = c\nu_{\bf v}(G) + d$.
For this reason, it may sometimes be convenient to assume that the smallest coordinate of ${\bf v}$ is $0$ and the largest coordinate is $1$
(or that ${\bf v}={\bf 1}$). Notice also that by dilation and translation with a constant vector, once can transform ${\bf v}$ to a nonnegative vector whose coordinate
sum is $1$, namely a probability distribution on $\F_k$.

Given a vector ${\bf v} \in {\mathbb R}^{|\F_k|}$, the extremal graph-theoretic question of interest is how small can $\nu_{\bf v}(G)$ be\footnote{
Maximizing $\nu_{\bf v}(G)$ is trivial. It is just the largest coordinate of ${\bf v}$, as if $K_n$ is $k$-decomposable, we can replace each copy of $K_k$ in
a $K_k$-decomposition of $K_n$ with a copy of $H$ where ${\bf v}_H$ is the maximum coordinate of ${\bf v}$ and the obtained graph $G$ has
$\nu_{\bf v}(G) = {\bf v}_H$.}. Thus, let $\nu_{\bf v}(n)$ denote the minimum of $\nu_{\bf v}(G)$ taken over all graphs $G$ with $n$ vertices
such that $K_n$ is $k$-decomposable. To formally extend this sequence to all $n$, one can synthetically define $\nu_{\bf v}(n) = \nu_{\bf v}(m)$
such that $m \le n$ is the largest integer such that $K_m$ is $k$-decomposable (trivially $K_1$ is $k$-decomposable).
It is not difficult to prove that the sequence $\nu_{\bf v}(n)$ converges (as shown later in this paper), but, as noted earlier, in most cases it is difficult,
and possibly intractable, to determine the limit. So, let $\nu_{\bf v} = \lim_{n \rightarrow \infty} \nu_{\bf v}(n)$.

To state our first result, we need to recall the notion of a fractional $k$-decomposition.
Let $\binom{G}{k}$ denote the set of $\binom{n}{k}$ $k$-vertex induced subgraphs of an $n$-vertex graph $G$.
For a pair of vertices $x,y$ of $G$, let $\binom{G}{k,x,y}$ be the $k$-vertex induced subgraphs of $G$ that contain both $x$ and $y$.
A {\em fractional $k$-decomposition} is a function $f:\binom{G}{k} \rightarrow [0,1]$ such that for each pair of vertices $x,y$, $\sum_{H \in \binom{G}{k,x,y}}f(H)=1$.
Clearly, a $k$-decomposition is also a fractional $k$-decomposition whose image is $\{0,1\}$.
Observe that every graph with $n \ge k$ vertices has a fractional $k$-decomposition, regardless of being $k$-divisible.
For ${\bf v} \in {\mathbb R}^{|\F_k|}$ indexed by $\F_k$
and a fractional $k$-decomposition $f$ of $G$, let $D_{\bf v}(f) = \sum_{H \in \F_k} {\bf v}_H f(H)$ where $f(H)$ is the sum of the values of $f$ on elements of
$\binom{G}{k}$ that are isomorphic to $H$. As before, it will be slightly more convenient to consider the normalized value
$\nu^*_{\bf v}(f) = D_{\bf v}(f)\binom{k}{2}/\binom{n}{2}$.
We therefore define the {\em optimal} fractional $k$-decomposition of $G$ with respect to ${\bf v}$ by $\nu^*_{\bf v}(G) = \max_{f} \nu^*_{\bf v}(f)$
where the maximum is taken over all fractional $k$-decompositions of $G$.
We define $\nu^*_{\bf v}(n)$ to be the minimum of $\nu^*_{\bf v}(G)$ taken over all graphs $G$ with $n$ vertices.
It is easy to verify that the sequence $\nu^*_{\bf v}(n)$ is non-decreasing and is upper bounded by the largest coordinate of ${\bf v}$, thus
let $\nu^*_{\bf v} = \lim_{n \rightarrow \infty} \nu^*_{\bf v}(n)$.
By the previous remark, we always have $\nu^*_{\bf v}(G) \ge \nu_{\bf v}(G)$, and consequently $\nu^*_{\bf v}(n) \ge \nu_{\bf v}(n)$ and
$\nu^*_{\bf v} \ge \nu_{\bf v}$.
The following is our first main result. We state it also for the analogous versions of tournaments and edge-colored graphs.
\begin{theorem}\label{t:main-1}
Let $k \ge 3$ be a given integer.
\begin{enumerate}
\item
Let ${\bf v} \in {\mathbb R}^{|\F_k|}$ be a given vector indexed by $\F_k$.
Then, $\nu_{\bf v} = \nu^*_{\bf v}$.
\item
Let ${\bf v} \in {\mathbb R}^{|\T_k|}$ be a given vector indexed by $\T_k$.
Then, $\nu_{\bf v} = \nu^*_{\bf v}$.
\item
Let $C$ be a finite set of colors and let ${\bf v} \in {\mathbb R}^{|\C_k|}$ be a given vector indexed by $\C_k$.
Then, $\nu_{\bf v} = \nu^*_{\bf v}$.
\end{enumerate}
In all cases, if $G$ has $n$ vertices such that $K_n$ is $k$-decomposable, then a $k$-decomposition $L$ of $G$ satisfying $\nu_{\bf v}(L) \ge \nu_{\bf v} - o_n(1)$ can
be constructed in polynomial time.
\end{theorem}
Note that the first case (that of $\F_k$) is equivalent to the special instance of the third case when the color set is $C=\{red,blue\}$.
Indeed, for a graph in $\F_k$ we can color its edges blue and its non-edges red thereby obtaining $\C_k$, and when considering an $n$-vertex  graph $G$ and an optimal $\F_k$-decomposition of it,
we can equivalently consider the optimal $\C_k$-decomposition of the blue-red coloring of $K_n$ where the edges of $G$ are colored blue and its non-edges are colored red.
It therefore suffices to prove only cases 2 and 3 of Theorem \ref{t:main-1}.

The proof of Theorem \ref{t:main-1} consists of two main ingredients. We first use a result from \cite{yuster-2005} which can, in particular, be formulated as follows. 
Given a family $\F$ of graphs, and given a fractional $\F$-decomposition of $G$ (assuming there is one), one can find a set $P$ of subgraphs of $G$ such that each element of $P$
is isomorphic to an element of $\F$ and any pair of vertices of $G$ is in at most one element of $P$ (if this pair is an edge of $G$, then it is also an edge
in the element of $P$ in which  it appears). Furthermore, the number of pairs that are not covered by $P$ is $o(|V(G)|^2)$.
So, assuming $G$ is dense, $P$ is a packing of elements of $\F$ in $G$ such that
almost all pairs of vertices of $G$ are packed and in this sense, it is an ``almost'' $\F$-decomposition.
The result in \cite{yuster-2005} extended an earlier result of Haxell and Ro\"dl \cite{HR-2001} where $\F$ is a single graph. Both results are actually more general,
as they show that any fractional packing (which may be far from a decomposition) can be converted to an integral packing with relatively small loss.
If we apply this result for $\F=\F_k$ we are close to proving the first part of Theorem \ref{t:main-1}, but there are two caveats.

Since we now have a vector ${\bf v} \in {\mathbb R}^{|\F_k|}$ associated, even if we start with an optimal fractional decomposition attaining $\nu^*_{\bf v}(G)$
it could be that the obtained integral ``almost decomposition'' distributes the weights to the elements of $\F_k$ in a way that decreases the total weight significantly below
$\nu^*_{\bf v}(G)$. However, fortunately, the proof from \cite{yuster-2005} (implicitly) shows that we can almost maintain the correct distribution.

The second (and more difficult) problem is that the obtained almost decomposition needs now to be modified to a full decomposition, and {\em without} affecting much the total weight, staying close to $\nu^*_{\bf v}(G)$. To this end, we use a fundamental result of Barber et al. \cite{BKLO-2016} based on the method of iterative absorption which enables us to achieve this goal, but with a price. To achieve their setting, we need, in fact, to first decompose a graph with high minimum degree to edge disjoint copies of some $K_m$ (here $m$ is huge compared to $k$, but fixed), and apply the aforementioned result of \cite{yuster-2005} to each element of this $K_m$-decomposition separately
(more precisely, to the subgraph of $G$ induced by the vertices of that element).
We then need to sparsify our obtained packing in order to achieve a setting suitable for the application of \cite{BKLO-2016}. 
The second and third part of Theorem \ref{t:main-1} are obtained using analogues of \cite{yuster-2005} for tournaments and edge-colored graphs.

Theorem \ref{t:main-1} provides a convenient mechanism to study certain natural decomposition problems, as it is sometimes much easier to obtain bounds
for the fractional problem. In fact, in many cases, we can glue optimal fractional decompositions of small graphs into good fractional decompositions for
arbitrary large graphs. We next give two very natural applications, but one may construct additional.
\begin{theorem}\label{t:applic-1}
Let $n \equiv 1,3 \bmod 6$. Any tournament on $n$ vertices has a triangle decomposition where the number of directed triangles in the decomposition is only
$0.0222n^2(1+o(1))$.
\end{theorem}
Let $\F^* \subset \F_k$. Consider the vector ${\bf v} \in {\mathbb R}^{|\F_k|}$ which assigns $0$ to the elements of $\F^*$ and $1$ to the elements of
$\F_k \setminus \F^*$.
We say that $\F^*$ is {\em essentially avoidable} if $\nu_{\bf v}=1$.
In other words, for every $k$-decomposable graph $G$, there is a $k$-decomposition of $G$ which almost completely avoids using elements from $\F^*$ (i.e. the fraction of elements of this
decomposition which are isomorphic to elements of $\F^*$ is $o_n(1))$. If $\F^*=\{H\}$, we say that $H$ is essentially avoidable.
A result from \cite{yuster-2019}, together with the proof of Theorem \ref{t:main-1} implies the following.
\begin{theorem}\label{t:applic-2}
$ $
\begin{enumerate}
\item
$C_5$ is essentially avoidable. More generally, if $k$ is odd and $\F^* \subset \F_k$ is the family of all graphs $H$ on $k$ vertices such that both $H$ and its complement are
Eulerian, then $\F^*$ is essentially avoidable.
\item
Almost all graphs are essentially avoidable. Namely, if ${\mathcal U}_k \subset \F_k$ is the set of all graphs on $k$ vertices that are not essentially avoidable, then
$|{\mathcal U}_k| = o(|\F_k|)$.
\end{enumerate}
\end{theorem}

For the first nontrivial case $k=3$, it is possible to determine $\nu_{\bf v}$ for every binary vector and certain additional types of vectors ${\bf v} \in {\mathbb R}^{|\F_3|}$ (these vectors are four dimensional as $|\F_3|=4$).
However, even for $k=3$, there are still some types of vectors for which we do not know $\nu_{\bf v}$.
For the case $k=4$ the situation is even more involved as we still do not know $\nu_{\bf v}$ even for all binary vectors.
We elaborate more on this in Section 4 which considers the small cases $k=3,4$.

It also seems plausible to try and evaluate the asymptotic behavior of $\nu_{\bf v}(G)$, namely, the asymptotic value of $\nu_{\bf v}(G)$ when $G \sim \G(n,p)$
is a random graph. In this case, it turns out that the problem can be completely solved, and the asymptotic value efficiently computed,
for all ${\bf v}$ and all constant $0 < p < 1$, as we prove in Section 5.

Our road-map follows. Section 2 contains the proof of Theorem \ref{t:main-1}. Our demonstrative applications, Theorem \ref{t:applic-1} and Theorem
\ref{t:applic-2} are the theme of Section 3. Section 4 focuses on small cases. Section 5 analyzes $\nu_{\bf v}(G)$ when $G \sim \G(n,p)$.

\section{Integer and fractional vector valued decompositions}

As noted in the introduction, it suffices to prove the second and third parts of Theorem \ref{t:main-1}. In this section we mostly prove the third part (the edge-coloring case).
The proof of the second part (the tournament case) follows along the same lines and requires only minor modifications, which are outlined in the last subsection of this section.

\subsection{From fractional decomposition to a similarly distributed integer packing}

Let $C$ be a finite set of colors and recall that $\C_k$ is the set of all edge colorings of $K_k$ with colors from $C$.
Suppose $G$ is a graph whose edges are colored by $C$. We call such $G$ a {\em $C$-colored graph} and
note that here we do not assume that $G$ is complete (so non-edges of $G$ correspond to non-colored pairs).
Let $\binom{G}{K_k}$ denote the set of $K_k$-subgraphs of $G$ and for $H \in \C_k$, let $\binom{G}{H} \subseteq \binom{G}{K_k}$ be the set of $K_k$-subgraphs of $G$ that
are color-isomorphic to $H$. More formally, for each $X \in \binom{G}{H}$ there is a bijection $b: V(X)\rightarrow V(H)$ such that $(u,v) \in E(X)$ and
$(b(u),b(v)) \in E(H)$ have the same color.

We can naturally extend the notion of a fractional $K_k$-decomposition to graphs that are not necessarily complete, such as $G$ above.
For an edge $e \in E(G)$, let $\binom{G}{K_k,e} \subseteq \binom{G}{K_k}$ be the set of $K_k$-subgraphs of $G$ that contain $e$.
We say that a function $f: \binom{G}{K_k} \rightarrow [0,1]$ is a fractional $K_k$-decomposition if
$\sum_{X \in \binom{G}{K_k,e}} f(X) = 1$ holds for each $e \in E(G)$. Notice that a necessary (though not sufficient) requirement for
$G$ to have a fractional $K_k$-decomposition is that each edge of $G$ belongs to at least one $K_k$-subgraph of $G$.

Now, suppose $f$ is a fractional $K_k$-decomposition of $G$ and that $G$ is $C$-colored.
For $H \in \C_k$, let $f(G,H) = \sum_{X \in \binom{G}{H}} f(X)$.
Since $f$ is a fractional $K_k$-decomposition, we have that
\begin{equation}\label{e:1}
\sum_{H \in \C_K} f(G,H) = \frac{|E(G)|}{\binom{k}{2}}\;.
\end{equation}
The following lemma follows implicitly from a generalization of the proof of the main result of \cite{yuster-2005}.
\begin{lemma}\label{l:1}
Let $C$ be a  finite set of colors, let $k \ge 3$ be an integer and let $\epsilon > 0$. There exists $n_0=n_0(k,C,\epsilon)$ such that the following holds.
Suppose $G$ is $C$-colored and has $n > n_0$ vertices.
Let $f$ be a fractional $K_k$-decomposition of $G$. Then for every $H \in \C_k$ there is a set $P_H$ of induced subgraphs of $G$ that are
color-isomorphic to $H$, such that $|P_H| \ge f(G,H) - \epsilon n^2$. Furthermore any two elements of $P=\cup_{H \in \C_k}P_H$ intersect in at most one vertex.
\end{lemma}
Since the main result of \cite{yuster-2005} is not proved for the edge-colored case (it is only for uncolored graphs) and since in any case the (rather short) proof there implies
Lemma \ref{l:1} only implicitly, we present the proof of Lemma \ref{l:1} in Subsection \ref{ss:l1}.

Notice that $P$ in Lemma \ref{l:1} is, in particular, a packing of $G$ with pairwise edge-disjoint copies of $K_k$.
As the following corollary shows, if we take an optimal fractional decomposition with respect to some ${\bf v} \in {\mathbb R}^{|\C_k|}$ and apply
Lemma \ref{l:1} to it, we obtain an integral packing of $G$ with elements of $\C_k$ that is close to an optimal $\C_k$-decomposition with respect to ${\bf v}$.
To be more formal, for ${\bf v} \in {\mathbb R}^{|\C_k|}$ indexed by $\C_k$
and a fractional $K_k$-decomposition $f$ of a $C$-colored graph $G$,
let $D_{\bf v}(f) = \sum_{H \in \C_k} {\bf v}_H f(G,H)$. As before, after normalizing we define
$\nu^*_{\bf v}(f) = D_{\bf v}(f)\binom{k}{2}/|E(G)|$
and define the {\em optimal} fractional $K_k$-decomposition of $G$ with respect to ${\bf v}$ by $\nu^*_{\bf v}(G) = \max_{f} \nu^*_{\bf v}(f)$
where the maximum is taken over all fractional $K_k$-decompositions of $G$. If $G$ has no fractional $K_k$-decomposition, then define $\nu^*_{\bf v}(G)=0$.
\begin{corollary}\label{c:1}
Let $C$ be a  finite set of colors, let $k \ge 3$ be an integer, let ${\bf v} \in {\mathbb R}^{|\C_k|}$, and let $\gamma > 0$. 
There exists $N_{\ref{c:1}}=N_{\ref{c:1}}(k,C,\gamma,{\bf v})$ such that the following holds for all $C$-colored graphs $G$ with $n > N_{\ref{c:1}}$ vertices
which have a fractional $K_k$-decomposition.
For every $H \in \C_k$ there is a set $P_H$ of induced subgraphs of $G$ that are color-isomorphic to $H$, such that any two elements of $P=\cup_{H \in \C_k}P_H$ intersect in  at most one vertex. Furthermore,
\begin{eqnarray*}
(a) & |P| \ge \frac{|E(G)|-\gamma n^2}{\binom{k}{2}}\;.\\
(b) & \sum_{H \in \C_k} {\bf v_H}|P_H| \ge \frac{|E(G)|}{\binom{k}{2}}\nu^*_{\bf v}(G)-\gamma n^2\;.
\end{eqnarray*}
\end{corollary}
\Proof
Let $s = \sum_{H \in \C_k} {\bf v}_H$.
Define $\epsilon=\gamma/(\binom{k}{2}|\C_k|)$ if $s < 1$ else set $\epsilon=\gamma/(\binom{k}{2}|\C_k|s)$. Let
$N_{\ref{c:1}}(k,C,\gamma,{\bf v})=n_0(k,C,\epsilon)$ where the latter is the constant from Lemma \ref{l:1}.
Let $G$ be a $C$-colored graph having $n > N_{\ref{c:1}}$ vertices. If $G$ has no fractional $K_k$-decomposition, then there is nothing to prove, 
so assume that $f$ is an optimal fractional $K_k$-decomposition of $G$ with respect to ${\bf v}$, thus $\nu^*_{\bf v}(f)=\nu^*_{\bf v}(G)$.
By Lemma \ref{l:1}, for every $H \in \C_k$ there is a set $P_H$ of induced subgraphs of $G$ that are color-isomorphic to $H$
such that $|P_H| \ge f(G,H) - \epsilon n^2$ and any two elements of $P=\cup_{H \in \C_k}P_H$ intersect in  at most one vertex.
Now,
\begin{eqnarray*}
\sum_{H \in \C_k} {\bf v_H}|P_H| & \ge & \sum_{H \in \C_k} {\bf v_H}\left(f(G,H) - \epsilon n^2\right)\\
& = & \left(\sum_{H \in \C_k} {\bf v_H}f(G,H)\right)  - \epsilon n^2 s\\
& = & \frac{|E(G)|}{\binom{k}{2}}\nu^*_{\bf v}(G)  - \epsilon n^2 s\\
& \ge & \frac{|E(G)|}{\binom{k}{2}}\nu^*_{\bf v}(G) - \gamma n^2
\end{eqnarray*}
which proves (b). To see (a) we just use (\ref{e:1}) and
\begin{eqnarray*}
|P| = \sum_{H \in \C_k} |P_H| \ge \left(\sum_{H \in \C_k} f(G,H)\right)  - \epsilon n^2 |\C_k| \ge \frac{|E(G)|-\gamma n^2}{\binom{k}{2}}\;.
\end{eqnarray*}
\qed

\subsection{From packing to decomposition}

The following lemma is a major ingredient of the proof of Theorem \ref{t:main-1}.
Recall that an {\em equitable partition} of a graph $G$ into $q$ parts $\P=\{W_1,\ldots,W_q\}$ is a partition of $V(G)$ such that $||W_i|-|W_j|| \le 1$ for all $1 \le i < j \le q$.
The {\em Tur\'an graph} with $q$ parts, denoted by $T(n,q)$  is the complete $q$-partite graph on $n$ vertices where the parts form an equitable partition.
\begin{lemma}\label{l:turan}
Let $k \ge 3$ be an integer. Then there exists $q_{\ref{l:turan}}(k)$ such that for all $q \ge q_{\ref{l:turan}}$ there exist
$N_{\ref{l:turan}}(q,k)$ and $\gamma=\gamma_{\ref{l:turan}}(q,k)$ such that the following holds for all $n > N_{\ref{l:turan}}$ for which $K_n$ is $k$-divisible.
Let $G$ be a complete graph on $n$ vertices,
let $\P$ be an equitable partition of $G$ into $q$ parts and let $G[\P]$ be the $T(n,q)$ spanning subgraph of $G$ formed by the parts of $\P$.
Suppose $P$ is a packing of $G[\P]$ with pairwise edge-disjoint copies of $K_k$ such that at most $\gamma n^2$ edges of $G[\P]$ are uncovered by elements of
$P$. Then, there is a sub-packing  $P' \subseteq P$ such that $|P|-|P'| \le 8\sqrt{\gamma} n^2$ and there is a $K_k$-decomposition of $G$ that contains $P'$.
\end{lemma}
The proof of Lemma \ref{l:turan} mainly follows from the proof of the main result of Barber et al. \cite{BKLO-2016}. We prove it in Subsection \ref{ss:l-turan}.
We will also need the following result which states that a graph with large enough minimum degree
has a fractional $K_m$-decomposition.
\begin{lemma}\label{l:2}\cite{BKLMO-2017,dukes-2012,yuster-2005-1}
For every integer $m \ge 3$, there exists $\alpha=\alpha(m) < 1$ such that every graph on $n$ vertices and minimum degree at least $\alpha n$ has a fractional $K_m$-decomposition. \qed
\end{lemma}
The first bound for $\alpha$ was given in \cite{yuster-2005-1} who proved that $\alpha \le 1-1/(9m^{10})$. This was later improved in \cite{dukes-2012}
to $1-2/(9m^2(m-1)^2)$ and in \cite{BKLMO-2017} to $1-1/(10^4m^{3/2})$.
It is worth noting that recently, an even stronger version of Lemma \ref{l:2} has been proved by Barber et al. \cite{BKLO-2016}.
In particular, they have proved that if $n$ is sufficiently large, and an $n$-vertex graph with minimum degree at least $\alpha n$ is $K_k$-divisible, then it has a $K_m$-decomposition (with roughly the same $\alpha$ as the one required for the fractional $K_m$-decomposition). However, using this stronger version for Lemma \ref{l:2}
will not make a difference in our arguments that follow. As mentioned above, we will, however, need to use the result from \cite{BKLO-2016} later in a subtler setting in order to prove Lemma \ref{l:turan}.

Next, we need the following simple lemma.
\begin{lemma}\label{l:3}
The sequence $\nu^*_{\bf v}(n)$ is non-decreasing and bounded from above, hence the limit $\nu^*_{\bf v}$ exists. In particular, for every for every $\epsilon > 0$, there exists $m$ such that $\nu^*_{\bf v}(m) \ge \nu^*_{\bf v} - \epsilon$.
\end{lemma}
\Proof
Let $s$ be the maximum coordinate of ${\bf v}$. Let $f$ be a fractional $K_k$-decomposition of $G=K_n$. Then, $f$ has $\nu^*_{\bf v}(f) \le s$.
So, the sequence $\nu^*_{\bf v}(n)$ is bounded from above by $s$. Next, we show that $\nu^*_{\bf v}(n) \ge \nu^*_{\bf v}(n-1)$.
Let $G$ be a $C$-colored complete graph on $n$ vertices. For each $v \in V(G)$, let $G_v$ be the induced subgraph of $G$ on  $V(G) \setminus v$ and
let $f_v$ be an optimal fractional $K_k$-decomposition of $G_v$ with respect to ${\bf v}$. So, by definition $\nu^*_{\bf v}(f_v) \ge \nu^*_{\bf v}(n-1)$.

Next, define a fractional $K_k$-decomposition $f$ of $G$ as follows. For each induced $k$-vertex subgraph $X$ of $G$, let
$$
f(X) = \frac{1}{n-2}\sum_{v \in V(G) \setminus V(X)} f_v(X)\;.
$$
It is easy to verify that the sum of the weights corresponding to each pair of vertices is precisely  $1$ so $f$ is indeed
a fractional $K_k$-decomposition of $G$ and that
$$
\nu^*_{\bf v}(f) = \frac{1}{n}\sum_{v \in V(G)} \nu^*_{\bf v}(f_v) \ge \nu^*_{\bf v}(n-1)\;.
$$
Hence, $\nu^*_{\bf v}(G) \ge \nu^*_{\bf v}(n-1)$ implying that $\nu^*_{\bf v}(n) \ge \nu^*_{\bf v}(n-1)$ and that the sequence is non-decreasing.
\qed

The following lemma immediately implies the third part of Theorem \ref{t:main-1}.
\begin{lemma}\label{l:main}
Let $C$ be a  finite set of colors, let $k \ge 3$ be an integer, let ${\bf v} \in {\mathbb R}^{|\C_k|}$ be indexed by $\C_k$, and let $\epsilon > 0$.
Then there exist $N_{\ref{l:main}}=N_{\ref{l:main}}(\epsilon,{\bf v})$ such that the following holds.
Let $G$ be a $C$-colored complete graph which is $k$-divisible and with $n > N_{\ref{l:main}}$ vertices.
Then, $G$ is $K_k$-decomposable and, furthermore $\nu_{\bf v}(G) \ge (\nu^*_{\bf v}-\epsilon)(1 - \epsilon)-\epsilon$.
\end{lemma}
\Proof
First notice that the lemma indeed implies the third part of Theorem \ref{t:main-1} since on the one hand we always have
$\nu_{\bf v}(n) \le \nu^*_{\bf v}(n) \le \nu^*_{\bf v}$ and on the other hand, the lemma shows that for every $\epsilon > 0$, if $n$ is sufficiently large, then
$\nu_{\bf v}(n) \ge (\nu^*_{\bf v}-\epsilon)(1 - \epsilon)-\epsilon$. Hence the limit $\nu_{\bf v}$ exists and equals $\nu^*_{\bf v}$.

We next establish some constants that are required for the proof and for the definition of $N_{\ref{l:main}}$.
Let $\epsilon > 0$ be given as in the statement of the lemma.
Let $m \ge k$ be the smallest integer such that  $\nu^*_{\bf v}(m) \ge \nu^*_{\bf v} - \epsilon/2$. Notice that $m$ exists by Lemma \ref{l:3}. 
Let $\alpha = \alpha(m)$ be the constant from Lemma \ref{l:2}.
Let $q=\lceil \max \{ 2/(1-\alpha)\,,\, 5/\epsilon, q_{\ref{l:turan}}(k) \} \rceil$.
Let $s$ be the maximum of $1$ and the maximum coordinate of ${\bf v}$.
Let $\gamma = \min\{\epsilon^2/(1024s^2k^4)\,,\,\gamma_{\ref{l:turan}}(q,k)\}$.
Let $N_{\ref{l:main}}=\max\{N_{\ref{c:1}}(k,C,\gamma,{\bf v})\,,\,N_{\ref{l:turan}}(q,k)\}$.
Let $n > N_{\ref{l:main}}$.

Let $G$ be a $C$-colored complete graph on $n$ vertices which is $k$-divisible.
Consider some arbitrary equitable partition $\P$ of $G$ into $q$ parts.
Let $G[\P]$ denote the spanning subgraph of $G$ consisting of all edges whose endpoints are in distinct parts.
Notice that $G[\P]$ is no longer complete, but we still view $G[\P]$ as a $C$-colored graph, where the edges of $G[\P]$ retain their colors.
Clearly, the minimum degree of $G[\P]$ satisfies $\delta(G[\P]) \ge \lfloor n-n/q \rfloor$ since in $G[\P]$ each vertex is adjacent to all other vertices but those in its part. By our choice of $q$ we have that $\delta(G[\P]) \ge \alpha n$.
Hence, by Lemma \ref{l:2}, $G[\P]$ has a fractional $K_m$-decomposition, call it $g$.

Recall that $\binom{G[\P]}{K_m}$ denotes the set of all $K_m$-subgraphs of $G[\P]$. So, $g: \binom{G[\P]}{K_m} \rightarrow [0,1]$ is such that for each edge of $G[\P]$,
the sum of the values of $g$ over all elements of $\binom{G[\P]}{K_m}$ that contain the edge is $1$.
We now define, for each $X \in \binom{G[\P]}{K_m}$, a fractional $K_k$-decomposition, denoted by $f_X$.
We take $f_X$ to be an optimal fractional $K_k$-decomposition of $X$ with respect to ${\bf v}$ (notice that $f_X$ exists since $X$ is a complete $C$-colored graph
and $|X|=m \ge k$). Thus, $\nu^*_{\bf v}(f_X)=\nu^*_{\bf v}(X)$.

We next define a fractional $K_k$-decomposition of $G[\P]$ denoted by $f$, as follows.
Let $Y$ be some $K_k$-subgraph of $G[\P]$. Let
\begin{equation}\label{e:2}
f(Y) = \sum_{X \in \binom{G[\P]}{K_m}, V(X) \supset V(Y)} f_X(Y)g(X)\;.
\end{equation}
Notice that $f$ is indeed a fractional $K_k$-decomposition of $G[\P]$ since $f_X$ is such for every $X \in \binom{G[\P]}{K_m}$ and since $g$ is a fractional $K_m$-decomposition of $G[\P]$.

We next estimate $D_{\bf v}(f) = \sum_{H \in \C_k} {\bf v}_H f(G[\P],H)$. By (\ref{e:2}) we have:
\begin{eqnarray*}
D_{\bf v}(f) & = & \sum_{X \in \binom{G[\P]}{K_m}} g(X)D_{\bf v}(f_X)\\
& = & \sum_{X \in \binom{G[\P]}{K_m}} g(X)\nu^*_{\bf v}(f_X)\frac{\binom{m}{2}}{\binom{k}{2}}\\
& = & \sum_{X \in \binom{G[\P]}{K_m}} g(X)\nu^*_{\bf v}(X)\frac{\binom{m}{2}}{\binom{k}{2}}\\
& \ge & \sum_{X \in \binom{G[\P]}{K_m}} g(X)\nu^*_{\bf v}(m)\frac{\binom{m}{2}}{\binom{k}{2}}\\
& = & \nu^*_{\bf v}(m)\frac{\binom{m}{2}}{\binom{k}{2}}\frac{|E(G[\P])|}{\binom{m}{2}}\\
& \ge & \left(\nu^*_{\bf v}-\frac{\epsilon}{2}\right)\frac{|E(G[\P])|}{\binom{k}{2}}\;.
\end{eqnarray*}
Since $\nu^*_{\bf v}(f)=D_{\bf v}(f)\binom{k}{2}/|E(G[\P])|$ we obtain from the last inequality that
$$
\nu^*_{\bf v}(G[\P]) \ge \nu^*_{\bf v}(f) \ge \nu^*_{\bf v}-\frac{\epsilon}{2}.
$$

We now apply Corollary \ref{c:1} to the graph $G[\P]$, which we can do since it has $n > N_{\ref{l:main}} \ge N_{\ref{c:1}}(k,C,\gamma,{\bf v})$ vertices
and since $G[\P]$ has a fractional $K_k$-decomposition.
By the corollary, we obtain that for every $H \in \C_k$ there is a set $P_H$ of induced subgraphs of $G[\P]$ that are color-isomorphic to $H$,
such that any two elements of $P=\cup_{H \in \C_k}P_H$ intersect in  at most one vertex. Furthermore,
$$
\sum_{H \in \C_k} {\bf v_H}|P_H| \ge \frac{|E(G[\P])|}{\binom{k}{2}}\nu^*_{\bf v}(G[\P])-\gamma n^2
$$
and
\begin{equation}\label{e:4}
|P| \ge \frac{|E(G[\P])|-\gamma n^2}{\binom{k}{2}}\;.
\end{equation}
But recall that $E(G[\P])$ consists of all $\binom{n}{2}$ edges of $G$ except those which have both of their endpoints in the same part of $\P$.
Thus, $|E(G[\P])| \ge \binom{n}{2}-q\binom{\lceil n/q \rceil}{2} \ge \binom{n}{2}-n^2/q$. Also, we have already proved that
$\nu^*_{\bf v}(G[\P]) \ge \nu^*_{\bf v}-\frac{\epsilon}{2}$. We therefore obtain using $q \ge 5/\epsilon$ that
\begin{eqnarray}
\sum_{H \in \C_k} {\bf v_H}|P_H| & \ge & \frac{\binom{n}{2}-n^2/q}{\binom{k}{2}}\left(\nu^*_{\bf v}-\frac{\epsilon}{2}\right)-\gamma n^2 \nonumber\\
& \ge & \frac{\binom{n}{2}}{\binom{k}{2}}\left(1-\frac{\epsilon}{2}\right) \left(\nu^*_{\bf v}-\frac{\epsilon}{2}\right)-\gamma n^2\;. \label{e:3}
\end{eqnarray}
Now, recall that each element $X \in P$ is also an induced $K_k$-subgraph of our complete graph $G$.
Let $G'$ denote the spanning subgraph of $G$ consisting of all edges that are not covered by elements of $P$.
Clearly, $G'$ is $K_k$-divisible since both $G$ and the complement of $G'$ (which is the edge-disjoint union of $K_k$'s) are $K_k$-divisible.
Now, suppose first that it was possible to find a $K_k$-decomposition of $G'$. Hence, in this case, there is a $K_k$-decomposition of $G$ that contains $P$.
We would therefore obtain from (\ref{e:3}) that
\begin{equation}\label{e:5}
\nu_{\bf v}(G) \ge \frac{\sum_{H \in \C_k} {\bf v_H}|P_H|}{\binom{n}{2}/\binom{k}{2}}
\ge \left(1-\frac{\epsilon}{2}\right) \left(\nu^*_{\bf v}-\frac{\epsilon}{2}\right)-\frac{\epsilon}{2}\;.
\end{equation}
Unfortunately, we have no guarantee that $G'$ has a $K_k$-decomposition.
Suppose, however, that it was possible to modify $P$ just a bit, say, by removing just a few of the elements of $P$ so that after this change, the
corresponding remainder graph $G'$ would  have a $K_k$-decomposition. Then, almost the same bound for $\nu_{\bf v}(G)$
would apply, assuming that $\sum_{H \in \F_k} {\bf v_H}|P_H|$ did not change much after the modification.
Fortunately, this is possible, as a consequence of Lemma \ref{l:turan}, as follows.
We can apply Lemma \ref{l:turan} since by (\ref{e:4}) $P$ covers all but at most $\gamma n^2$ edges of $G[\P]$.
The lemma shows that there is a sub-packing $P' \subseteq P$ such that $|P|-|P'| \le 8\sqrt{\gamma} n^2$ and there is a $K_k$-decomposition $P^*$ of $G$ that contains $P'$.
Recall that $s$ is the maximum of $1$ and the maximum coordinate of ${\bf v}$. We therefore have by (\ref{e:5}) that:
$$
\nu_{\bf v}(G) \ge \left(1-\frac{\epsilon}{2}\right) \left(\nu^*_{\bf v}-\frac{\epsilon}{2}\right)-\frac{\epsilon}{2}- \frac{8s\sqrt{\gamma} n^2}{\binom{n}{2}/\binom{k}{2}}
\ge (1-\epsilon)(\nu^*_{\bf v}-\epsilon)-\epsilon
$$
where we have used that $\sqrt{\gamma} \le  \epsilon/(32sk^2)$.
\qed

Lemma \ref{l:main} can be implemented in polynomial time as claimed in the statement of Theorem \ref{t:main-1}. Namely the $K_k$-decomposition $P^*$ in the lemma can be constructed in time which is polynomial in $n=|V(G)|$. To see this, we first observe that Lemma \ref{l:1} can be implemented in polynomial time (i.e. constructing the packing $P$ in that lemma), as proved in \cite{yuster-2005}. This implies that Corollary \ref{c:1} can be implemented in polynomial time, since finding the 
optimal fractional $K_k$-decomposition of $G$ with respect to ${\bf v}$ denoted by $f$ in the proof of the corollary can be found in polynomial time
using linear programming (the number of variables is $O(n^k)$ as the number of $K_k$ and the number of constraints is only $O(n^2)$ as the number of edges).
Once we obtain the packing $P$ of Corollary \ref{c:1}, we apply Lemma \ref{l:turan} which constructs $P^*$ in polynomial time, as Lemmas \ref{l:absorber} and
\ref{l:almost-decomp} in Subsection \ref{ss:l-turan} can be implemented in polynomial time as proved in \cite{BKLO-2016}.

\subsection{Proof of Lemma \ref{l:turan}}\label{ss:l-turan}

The proof of Lemma \ref{l:turan} is based on the proof of the main result of \cite{BKLO-2016} (Theorem 1.3 there).
In fact, we will only need to use a special case of that result, for the case of the small graph being $K_k$ and for the case of the host graph being $G=K_n$
although most of the arguments in \cite{BKLO-2016} are still required even for this special case, which is not surprising 
since this special case implies Wilson's decomposition theorem for the case of $K_k$.
To achieve the setting in \cite{BKLO-2016} we require some definitions taken from there.

For a graph $G$, a positive integer $q$ and a real $\delta > 0$,
a {\em $(q,\delta)$-partition} of $G$ is an equitable partition
$\P = \{V_1 \ldots,V_q\}$ of $V(G)$ such that for each $1 \le i \le q$
and each $v \in V (G)$, $d_G(v,V_i) \ge \delta|V_i|$. Here $d_G(v,V_i)$ denotes the
number of neighbors of $v$ in $V_i$. Notice that if $G=K_n$, then $G$ trivially has a $(q,\delta)$-partition, but a straightforward probabilistic argument shows that
this also holds if $G$ is just an $n$-vertex graph with minimum degree slightly larger than $\delta n$ and $n$ is sufficiently large (Proposition 7.3 in \cite{BKLO-2016}).
For an equitable partition $\P$ into $q$ parts, recall that $G[\P]$ denotes the $q$-partite subgraph of $G$ induced by the parts of $\P$.

For an equitable partition $\P$ into $q$ parts, a {\em refinement} of $\P$ is obtained by taking an equitable partition into $q$ parts of each part of $\P$.
Notice that a refinement is an equitable partition into $q^2$ parts.

Let $\P_1$ be an equitable partition of $V(G)$ and for each $2 \le i \le \ell$ let $\P_i$ be a refinement of $\P_{i-1}$.
We say that $\P_1,\ldots,\P_\ell$ is a {\em $(q,\delta,m)$-partition sequence} of $G$ if the following hold.\\
(i) $\P_1$ is a $(q,\delta)$-partition of $G$.\\
(ii) For each $2 \le i \le \ell$ and each $V \in \P_{i-1}$, $\P_i[V]$ is a $(q,\delta)$-partition of $G[V]$.\\
(iii) Each part of $\P_\ell$ is of size $m$ or $m-1$.\\
Once again, if $G=K_n$, then it trivially has a $(q,\delta,m)$-partition sequence, but also if $n$ is sufficiently large
it is easy to prove that a graph with $n$ vertices and minimum degree slightly larger than $\delta n$ (say minimum degree at least $(\delta+\epsilon)n$)
has a $(q,\delta,m)$-partition sequence, where $m$ is bounded by a constant depending only on $q$ and $\epsilon$ (Lemma 7.4 in \cite{BKLO-2016}).

The first major ingredient in the proof of Theorem 1.3 in \cite{BKLO-2016} is that of the existence of an {\em absorber}.
Informally, an absorber $A^*$ of a $K_k$-divisible graph $G$ is a $K_k$-divisible spanning subgraph of $G$ with small maximum degree which has the following property.
Suppose we take an ``almost $K_k$-decomposition'' of the spanning subgraph $G'$ obtained from $G$ after removing the edges of $A^*$.
Let $H^*$ be the leftover edges of $G'$ uncovered by the almost decomposition. Note that $H^*$ is also $K_k$-divisible.
Then $A^*$ has the property that $A^* \cup H^*$ has a $K_k$-decomposition (and hence so does $G$).
Of course, in order to obtain such an $A^*$ we need to make sure that the set of possible $H^*$ is small (in particular, if one can guarantee that $H^*$ has 
no more than $O(n)$ edges, this will limit the number of possibilities for $H^*$).
The formal definition of such an absorber is given in Lemma 8.1 there, which is stated here for the special case of $K_k$.
Note that some of the notations have been changed to adjust to the notations in the present paper.

\begin{lemma}\label{l:absorber}
Let $k \ge 3$ and $\epsilon > 0$. Then there exists $m_{\ref{l:absorber}}(\epsilon,k)$ such that the following holds for all $m \ge m_{\ref{l:absorber}}$.
There exists  $N_{\ref{l:absorber}}(\epsilon,k,m)$ such that for all $n > N_{\ref{l:absorber}}$ the following holds.
Set $\delta \coloneqq 1-1/(3k)+\epsilon$, $t \coloneqq \lceil n/m \rceil$ and let $G$ be a graph with $n$ vertices.
Let $\P=\{V_1,\ldots,V_t\}$ be an equitable partition of $V(G)$ so that each part has size $m$ or $m-1$. Suppose that $\delta(G[\P]) \ge \delta n$
and $\delta(G[V_i]) \ge \delta|V_i|$ for each $1 \le i \le t$. Then $G$ contain a $K_k$-divisible subgraph $A^*$ such that:\\
(i) $\Delta(A^*[\P]) \le \epsilon^2 n$ and $\Delta(A^*[V_i]) < k$ for each $1 \le i \le t$.\\
(ii) If $H^*$ is a $K_k$-divisible graph on $V(G)$ that is edge-disjoint from $A^*$ and has $E(H^*[\P]) = \emptyset$, then $H^* \cup A^*$ has a $K_k$-decomposition. \qed
\end{lemma}

For a subgraph $X$ of $G$ let $G-X$ denote the spanning subgraph of $G$ obtained by removing the edges of $X$.
In order to apply Lemma \ref{l:absorber} one first needs to decompose $G - (A^* \cup H^*)$. This is the other major ingredient in \cite{BKLO-2016},
which appears as Lemma 10.1 there. The following is a version of Lemma 10.1 for the special case of $K_k$ and with an addendum that follows from its proof.
\begin{lemma}\label{l:almost-decomp}
Let $k \ge 3$ and $\epsilon > 0$. Then there exists $q_{\ref{l:almost-decomp}}(k,\epsilon)$ such that the following holds for all $q \ge q_{\ref{l:almost-decomp}}$.
There exists $\gamma_{\ref{l:almost-decomp}}(q,\epsilon,k)$ such that the following holds for all $\gamma \le \gamma_{\ref{l:almost-decomp}}$, for all
$m \ge m_{\ref{l:almost-decomp}}(\gamma)$ and for every $K_k$-divisible graph
$G$ on $n$ vertices. Define $\delta \coloneqq \max\{\alpha(k),1-1/(3k)\}$ where $\alpha(k)$ is the constant from Lemma \ref{l:2}. Suppose $\P_1,\ldots,\P_\ell$ is a
$(q,\delta+\epsilon,m)$-partition sequence of $G$. Then there exists a subgraph $H^*$ of $\cup_{V \in \P_\ell} G[V]$ such
that $G - H^*$ has a $K_k$-decomposition $P^*$. Furthermore, if $P$ is packing of $G[\P_1]$ covering all but at most $2\gamma n^2$ edges of $G[\P_1]$,
then there exists such a $P^*$ such that $|P^* \setminus P| \le 6\sqrt{\gamma} n^2$.
\qed
\end{lemma}
We note that the ``Furthermore'' part does not appear in the statement of Lemma 10.1 in \cite{BKLO-2016}, but immediately follows from its proof.
Indeed, the first part of the proof (Lemma 10.6 there) proceeds as follows. Take any packing $P$ of $G[\P_1]$ that covers all but at most $2\gamma n^2$ edges of $G$.
By removing at most $6\sqrt{\gamma} n^2$ elements from $P$ you obtain a packing $P'$ such that the subgraph $H$ of $G$ consisting of the edges of $G$ that are uncovered by
$P'$ has some nice properties (stated as (G1) and (G2) in Lemma 10.6). From there onwards the proof Lemma 10.1 proceeds by iteratively improving $P'$ in $\ell-1$ steps
where in step $i$ one obtains an almost optimal packing in each part of $\P_i$ which covers also the remaining uncovered edges between parts of the previous
partition $\P_{i-1}$ until obtaining $P^*$ of Lemma \ref{l:almost-decomp}. In particular, $P^*$ still retains almost all of the element of the initial packing $P$,
but at most $6\sqrt{\gamma} n^2$ elements.

\vspace{4pt}
\noindent
{\bf Proof of Lemma \ref{l:turan}:}
Let $k \ge 3$ be an integer. Define the following constants.\\
(i) $\delta \coloneqq \max\{\alpha(k)\,,\,1-1/(3k)\}$ where $\alpha(k)$ is the constant from Lemma \ref{l:2}.\\
(ii) $\epsilon = (1-\delta)/10$ and $\epsilon'=\epsilon/3$.\\
(iii) $q_{\ref{l:turan}}(k)= \max \{ 2/(1-\delta-\epsilon')\,,\, q_{\ref{l:almost-decomp}}(k,\epsilon') \} $ and let $q \ge q_{\ref{l:turan}}$.\\
(iv) $\gamma= \gamma_{\ref{l:turan}}(q,k)=\gamma_{\ref{l:almost-decomp}}(q,\epsilon',k)$.\\
(v) $m=\max\{ m_{\ref{l:almost-decomp}}(\gamma)\,,\,m_{\ref{l:absorber}}(\gamma,k)\}$.\\
(vi) $N_{\ref{l:turan}}(q,k) = N_{\ref{l:absorber}}(\gamma,k,m)$.

Now let $n > N_{\ref{l:turan}}$ such that $K_n$ is $k$-divisible.
Let $G$ be a complete graph on $n$ vertices,
let $\P$ be an equitable partition of $G$ into $q$ parts and let $G[\P]$ be the $T(n,q)$ spanning subgraph of $G$ formed by the parts of $\P$.
Suppose $P$ is a packing of $G[\P]$ with pairwise edge-disjoint copies of $K_k$ such that at most $\gamma n^2$ edges of $G[\P]$ are uncovered by elements of
$P$.

Suppose $\P_1,\ldots,\P_\ell$ is a
$(q,\delta+\epsilon,m)$-partition sequence of $G$ where $\P_1=\P$. Observe that such a $(q,\delta+\epsilon,m)$-partition exists since $G$ is a complete graph, since
$\P_1=\P$ is an equitable partition into $q$ parts, and since $\delta+\epsilon < 1$.

Let $G_1=G[\P]=G[\P_1]$ and let $G_{\ell+1}=G-G[P_\ell]$. So, $G_{\ell+1}$ consists of all edges with both endpoints in the same part of $\P_\ell$. 
Consider now the graph $H=G_1 \cup G_{\ell+1}$ (i.e. the spanning subgraph of $G$ consisting of all the edges of $G_1$ and $G_{\ell+1}$)
and consider the partition $\P_\ell$ of $H$. First observe that the minimum degree $\delta(H[\P_\ell])$ is at least
$n-\lceil n/q \rceil \ge (\delta+\epsilon')n$ where we have used here that $q \ge  2/(1-\delta-\epsilon')$.
Similarly, for each $V \in \P_\ell$ we have $\delta(H[V]) =|V|-1 \in \{m-1,m-2\}$. So, $\delta(H[V]) \ge (\delta+\epsilon') |V|$
since $\delta+\epsilon' < 1$. We may therefore apply Lemma \ref{l:absorber} where $H$ plays the role of $G$, $\gamma$ plays the role of $\epsilon$
and $\P_\ell$ plays the role of $\P$.

By Lemma \ref{l:absorber}, $H$ contains a $K_k$-divisible subgraph $A^*$ such that:\\
(i) $\Delta(A^*[\P_\ell]) \le \gamma^2 n$ and $\Delta(A^*[V]) < k$ for each $V \in \P_\ell$.\\
(ii) If $H^*$ is a $K_k$-divisible graph on $V(G)=V(H)$ that is edge-disjoint from $A^*$ and $E(H^*[\P_\ell]) = \emptyset$, then $H^* \cup A^*$ has a $K_k$-decomposition.\\
Observe that (i) and (ii) imply also that $\Delta(A^*) < \gamma^2 n+k$.
Let $G'=G-A^*$. Thus, $G'$ is also $K_k$-divisible.
Note that for each $V \in \P_1$ and each $v \in V(G)$ we have $d_{G'}(v,V) \ge d_{G}(v,V)-\Delta(A^*) \ge (|V|-1)-(\gamma^2 n+k-1) \ge (\delta+\epsilon')|V|$.
So, $\P_1$ is a $(q,\delta+\epsilon')$-partition of $G'$.
Note also that by (i) we have that $\Delta(A^*-A^*[\P_1]) < k$ so $\P_1,\ldots,\P_\ell$ is also a $(q,\delta+\epsilon',m)$-partition sequence of $G'$.
Recall also that the packing $P$ covered at most $\gamma n^2$ edges of $G[\P]$. Let $P'' \subset P$ be the elements of $P$ which are entirely in $G'$.
Hence, each element of $P \setminus P''$ contain an edge of $A^*$. Since $\Delta(A^*) < \gamma^2 n+k$, we have that the number of edges of $A^*$ is at most
$\gamma^2n^2+nk$. It follows that $P''$ covers all elements of $G'(\P)$ but at most $\gamma n^2 + \binom{k}{2}(\gamma^2n^2+nk) < 2\gamma n^2$.

We can therefore apply Lemma \ref{l:almost-decomp} to $G'$ playing the role of $G$, $\epsilon'$ playing the role of $\epsilon$, and
$P''$ playing the role of $P$ in that lemma. By Lemma \ref{l:almost-decomp} we obtain a
subgraph $H^*$ of $\cup_{V \in \P_\ell} G'[V]$ such
that $G' - H^*$ has a $K_k$-decomposition $P^*$. Furthermore, $|P^* \setminus P''| \le 6\sqrt{\gamma} n^2$.
But now, by (ii) $A^* \cup H^*$ has a $K_k$-decomposition, so together with $P^*$ this forms a $K_k$-decomposition of $G=K_n$ containing all
but at most $6\sqrt{\gamma} n^2$ elements of $P''$ thus all but at most $6\sqrt{\gamma} n^2+2\gamma n^2 \le 8\sqrt{\gamma}n^2$ elements of $P$.
\qed

\subsection{Proof of Lemma \ref{l:1}}\label{ss:l1} 

As noted earlier, Lemma \ref{l:1} follows implicitly from the main result in \cite{yuster-2005}.
That result is stated in terms of uncolored graphs, while here we need the colored version.
Thus, we reproduce the arguments in the proof of \cite{yuster-2005} where the lemmas there whose proofs remain identical or
for which the colored version is an immediate extension are only restated in their colored version without proof,
but with reference to the original lemma in \cite{yuster-2005}.

We first need to recall the edge-colored version of the Szemer\'edi's regularity
lemma \cite{szemeredi-1978}. Let $G=(V,E)$ be a $C$-colored graph, and let $A$ and $B$ be two disjoint subsets
of $V(G)$. If $A$ and $B$ are non-empty and $c \in C$, let $E_c(A,B)$ denote the set of edges between them that are colored $c$.
The {\em $c$-density} between $A$ and $B$ is defined as
$$
d_c(A,B) = \frac{|E_c(A,B)|}{|A||B|}.
$$
For $\gamma>0$ the pair $(A,B)$ is called {\em $\gamma$-regular}
if for every $X \subseteq A$ and $Y \subseteq B$ satisfying
$|X|\ge \gamma |A|$ and $|Y| \ge \gamma |B|$ we have
$$
|d_c(X,Y)-d_c(A,B)| \le \gamma ~~~{\rm for~all~} c \in C\;.
$$
An equitable partition of the set of vertices $V$ of a $C$-colored graph $G$ into the classes $V_1,\ldots,V_m$ is
called {\em $\gamma$-regular} if all but at most $\gamma \binom{m}{2}$ of the pairs
$(V_i,V_j)$ are $\gamma$-regular. The regularity lemma (colored version) states the following:
\begin{lemma}
\label{l:21}
Let $C$ be a finite set of colors and let $\gamma>0$. There is an integer $M(\gamma,C)>0$ such that for
every $C$-colored graph $G$ of order $n > M$ there is a $\gamma$-regular partition of
the vertex set of $G$ into $m$ classes, for some $1/\gamma < m < M$. \qed
\end{lemma}
The proof of Lemma \ref{l:21} is completely analogous to the proof of the original regularity lemma.

For an edge $(x,y)$ of a $C$-colored graph, let $c(x,y)$ denote its color.
Let $H$ be a $C$-colored graph with $V(H)=\{1,\ldots,k\}$, $k \geq 3$.
Let $W$ be a $C$-colored $k$-partite graph with vertex classes $V_1,\ldots,V_k$. A subgraph
$J$ of $W$ with $V(J)=\{v_1,\ldots,v_k\}$ is {\em partite-color-isomorphic} to $H$ if $v_i \in V_i$
for $i=1,\ldots,k$ and the map $i \rightarrow v_i$ is a {\em color preserving isomorphism} from $H$ to $J$.
Namely, $(i,j) \in E(H)$ if and only if $(v_i,v_j) \in E(J)$ and in case they are both edges, then $c(i,j)=c(v_i,v_j)$.

The following is a standard counting lemma whose proof follows from the definition of $\gamma$-regularity.
It is analogous to Lemma 2.2 of \cite{yuster-2005}.
\begin{lemma}
\label{l:22}
Let $C$ be a finite set of colors, let $k \ge 3$ be a positive integer, and let $\delta$ and $\zeta$ be positive reals.
There exist $\gamma=\gamma(\delta,\zeta,k,C)$ and $T=T(\delta,\zeta,k,C)$ such that the following holds.
Let $H$ be a $C$-colored graph with $V(H)=\{1,\ldots,k\}$
and let $W$ be a $C$-colored $k$-partite graph with vertex classes $V_1,\ldots,V_k$ where $|V_i|=t > T$ for $i=1,\ldots,k$.
Furthermore, for each $(i,j) \in E(H)$, $(V_i,V_j)$ is
a $\gamma$-regular pair with $d_{c(i,j)}(V_i,V_j) \geq \delta$ and for each $(i,j) \notin E(H)$,
$E(V_i,V_j)=\emptyset$.
Then, there exists a spanning subgraph $W'$ of $W$, consisting of at least $(1-\zeta)|E(W)|$ edges
such that the following holds.
For an edge $e \in E(W')$, let $count(e)$ denote the number of subgraphs of $W'$ that are partite-color-isomorphic
to $H$ and that contain $e$.
Then, for all $e \in E(W')$, if $e \in E(V_i,V_j)$, then
$$
\left|count(e) - t^{k-2} \frac{\prod_{(s,p) \in E(H)}d_{c(s,p)}(V_s,V_p)}{d_{c(i,j)}(V_i,V_j)}\right| < \zeta t^{k-2}.
$$
\qed
\end{lemma}

We need the result of Frankl and R\"odl \cite{FR-1985} on near perfect matchings of uniform hypergraphs.
Recall that if $x,y$ are two vertices of a hypergraph then $deg(x)$ denotes the degree of $x$ and $deg(x,y)$ denotes the number of edges that contain
both $x$ and $y$. We use the version of the Frankl and R\"odl Theorem due to Pippenger.
\begin{lemma}
\label{l:23}
For an integer $r \geq 2$ and a real $\beta > 0$ there exists $\mu=\mu(r,\beta) > 0$ so that:
If the $r$-uniform hypergraph $L$ on $q$ vertices has the following properties for some $d$:\\
(i) $(1-\mu)d < deg(x) < (1+\mu)d$ holds for all vertices,\\
(ii) $deg(x,y) < \mu d$ for all distinct $x$ and $y$,\\
then $L$ has a matching of size at least $(q/r)(1-\beta)$. \qed
\end{lemma}

Let $C$ be a  finite set of colors, let $k \ge 3$ be an integer and let $\epsilon > 0$.
Let $\delta=\beta=\epsilon/4$. Let $\mu=\mu(\binom{k}{2},\beta)$ be as in Lemma \ref{l:23}.
Let $\zeta=\mu \delta^{k^2}/2$. Let $\gamma=\gamma(\delta, \zeta, k,C)$ and
$T=T(\delta, \zeta, k,C)$ be as in Lemma \ref{l:22}.
Let $M=M(\gamma\epsilon/(25k^2),C)$ be as in Lemma \ref{l:21}.
Finally, we shall define $n_0=n_0(k,C,\epsilon)$ to be a sufficiently large constant, depending on the above chosen parameters,
and for which the inequalities stated in the proof below hold.

Fix an $n$-vertex $C$-colored graph $G$ with $n > n_0$ vertices and assume that $G$ has a fractional $K_k$-decomposition
$f:\binom{G}{K_k} \rightarrow [0,1]$.
We apply Lemma \ref{l:21} to $G$ and obtain a $\gamma'$-regular
partition with $m'$ parts, where $\gamma'=\gamma \epsilon/(25k^2)$ and $1/\gamma' < m' < M$.
Denote the parts by $U_1,\ldots,U_{m'}$.
Notice that the size of each part is either $\lfloor n/{m'} \rfloor$ or $\lceil n/{m'} \rceil$.
For simplicity we may and will assume that $n/{m'}$ is an integer, as this assumption does not
affect the asymptotic nature of the result. Similarly, we assume that
$25k^2/\epsilon$ and $n/(25m'k^2/\epsilon)$ are integers.

We randomly partition each $U_i$ into $25k^2/\epsilon$ equal parts of size $n/(25m'k^2/\epsilon)$ each.
All $m'$ partitions are independent. We now have $m=25m'k^2/\epsilon$ {\em refined} vertex classes, denoted
$V_1,\ldots,V_m$.
Suppose $V_i \subset U_s$ and $V_j \subset U_t$ where $s \neq t$. We claim that if
$(U_s,U_t)$ is a $\gamma'$-regular pair, then $(V_i,V_j)$ is a $\gamma$-regular pair.
Indeed, if $X \subseteq V_i$ and $Y \subseteq V_j$ have $|X|, |Y| \ge \gamma n/(25m'k^2/\epsilon)$, then
$|X|, |Y| \ge \gamma' n/m'$ and so $|d_c(X,Y) - d_c(U_s,U_t)| \le \gamma'$ for each $c \in C$. Also
$|d_c(V_i,V_j) - d_c(U_s,U_t)| \le \gamma'$. Thus, $|d_c(X,Y) - d_c(V_i,V_j)| \le 2\gamma' \le \gamma$.

Let $X$ be some $K_k$-subgraph of $G$. We call $X$ {\em good} if its $k$ vertices belong to distinct vertex classes of the refined partition.
Since the probability that two vertices of $X$ belong to the same vertex class of the refined partition is less than $\epsilon/(25k^2)$,
the probability that $X$ is not good is at most
$\binom{k}{2}\epsilon/(25k^2) < \epsilon/50$.
Since $f$ is a fractional $K_k$-decomposition, the sum of its values is $|f|=|E(G)|/\binom{k}{2} < n^2$.
Hence, if $f^{**}$ is the restriction of $f$ to
good elements (the non-good elements having $f^{**}(X)=0$), then the expected sum of the values of $f^{**}$ 
is at least $|f|(1-\epsilon/50)$. We therefore {\em fix} a partition $V_1,\ldots,V_m$ for which
$|f^{**}| \ge |f|(1-\epsilon/50)$. Notice that $f^{**}$ is no longer a fractional $k$-decomposition; it is merely
a fractional $K_k$-packing of $G$ (i.e. for each edge of $G$, the sum of the values of $f^{**}$ on the elements of $\binom{G}{K_k}$ that
contain the edge is at most $1$).
Furthermore, for each $H \in \C_k$ we have that
$$
f^{**}(G,H) \ge f(G,H) - (|f|-|f^{**}|) \ge f(G,H)-\frac{\epsilon}{50}|f| \ge f(G,H) - \frac{\epsilon}{50} n^2\;.
$$

Let $G^*$ be the spanning subgraph of $G$ consisting of the following edges:
An edge $(u,v) \in E(G)$ is in $E(G^*)$ if and only if $u \in V_i$, $v \in V_j$, $i \neq j$, $(V_i,V_j)$ is a $\gamma$-regular
pair, and $d_{c(u,v)}(V_i,V_j) \ge \delta$.
(thus, we discard edges inside classes, between non regular pairs, or if the color of the edge is sparse in the pair to which it belongs).
Let $f^*$ be the restriction of $f^{**}$ to copies of $K_k$ in $G^*$.
We claim that
$|f^*| > |f^{**}|-0.6\delta n^2$.
Indeed, by considering the number of discarded edges we get (using $\delta \gg \gamma' \ge 1/m'$)
\begin{eqnarray*}
|f^{**}| - |f^*| & \le & |E(G) - E(G^*)| \\
& < & \gamma' {\binom{m'}{2}}\frac{n^2}{{m'}^2} + {\binom{m'}{2}}(\delta+\gamma')\frac{n^2}{{m'}^2} + {m'}{\binom{n/{m'}}{2}}\\
& < & 0.6 \delta n^2\;.
\end{eqnarray*}
In particular, for each $H \in \C_k$ we have that
$$
f^{*}(G^*,H) \ge f^{**}(G,H) - (|f^{**}| - |f^*|) \ge f(G,H) - \frac{\epsilon}{50} n^2 - 0.6\delta n^2 \ge f(G,H) -\frac{\epsilon}{10} n^2\;.
$$

Let $R$ denote the $m$-vertex {\em multigraph} whose vertices are $\{1,\ldots,m\}$ and a pair $(i,j)$ with color $c$ is an edge of $R$
if and only if $(V_i,V_j)$ is a $\gamma$-regular pair and $d_c(i,j) \ge \delta$.
Notice that $R$ is indeed a multigraph but any two multiple edges have distinct colors. We define a
fractional $K_k$-packing $f'$ of $R$ as follows. Let $X$ be a
subgraph of $R$ that is color-isomorphic to some $H \in \C_k$ and assume that the vertices of $X$ are
$\{u_1,\ldots,u_k\}$ where $u_i$ plays the role of vertex $i$ in $H$.
We define $f'(X)$ to be the sum of the values of $f^*$ taken over all
subgraphs of $G^*[V_{u_1},\ldots,V_{u_k}]$ which are partite-color-isomorphic to $H$,
divided by $n^2/m^2$ (and where the isomorphism is $i \rightarrow u_i$). Notice that $|f'|=m^2 |f^*|/n^2$ since every $K_k$-subgraph of $G^*$ contributes its weight
(divided by $n^2/m^2$) to the sum of the weights of $f'$.
Likewise
$$
f'(R,H) = f^{*}(G^*,H)\frac{m^2}{n^2}\;.
$$

We use $f'$ to define a random partition of $E(G^*)$. Our parts correspond to the
copies of elements of $\binom{R}{K_k}$.
We denote the partition by $\Q=\{Q_X\,:\,X \in \binom{R}{K_k}\}$.
Let $X \in \binom{R}{K_k}$
and assume that $X$ contains the edge $(i,j)$ of $E(R)$ and that the color of the edge is $c$.
Each $e \in E_c(V_i,V_j)$ (which, by the definition of $R$, must be an edge of $G^*$) is chosen to
be in $Q_X$ with probability $f'(X)/d_c(V_i,V_j)$.
The choices made by distinct edges of $G^*$ are independent.
Notice that this random coloring is legal (in the sense that the sum of probabilities is
at most one) since the sum of $f'(X)$ taken
over all possible $X$ containing the edge $(i,j)$ of $E(R)$ whose color is $c$ is at most $d_c(V_i,V_j)$.
Notice also that some edges of $G^*$ might stay unassigned to a part in our random partitioning (as maybe
an edge $(i,j)$ of $E(R)$ whose color is $c$ does not belong to any $X$). In this case, we can assign such unassigned edges of $G^*$
to some ``spare part'', denoted $Q_0$, so that $\Q=\{Q_X\,:\,X \in \binom{R}{K_k}\} \cup \{Q_0\}$ is indeed a partition of $E(G^*)$.

Let $X$ be a
subgraph of $R$ that is color-isomorphic to some $H \in \C_k$, and assume that $f'(X) > m^{1-k}$ (we need this assumption in the lemmas below).
Without loss of generality, assume that the vertices of $X$ are $\{1,\ldots,k\}$ where $i \in V(X)$ plays
the role of $i \in V(H)$.
Let $W_X=G^*[V_1,\ldots,V_k]$.
Notice that $W_X$ is a subgraph of $G^*$ which satisfies the conditions in Lemma \ref{l:22},
since $t=n/m > n_0\epsilon/(25k^2M) > T$ (here we assume $n_0 > 25k^2MT/\epsilon$).
Let $W'_X$ be the spanning subgraph of $W_X$ whose existence is guaranteed in Lemma \ref{l:22}.
Let $Z_X$ denote the spanning subgraph of $W'_X$ consisting only of the edges that belong to the part $Q_X$.
Notice that $Z_X$ is a random subgraph of $W'_X$.
For an edge $e \in E(Z_X)$, let $S_X(e)$ denote the set of subgraphs of $Z_X$ that contain $e$
and that are partite-color-isomorphic to $H$. Put $s_X(e)=|S_X(e)|$. the proof of the following two lemmas are identical to the proofs of Lemmas 3.1 and Lemma 3.2
in \cite{yuster-2005}, respectively.
\begin{lemma} \label{l:31}
With probability at least $1-m^3/n$, for all $e \in E(Z_X)$,
$$
\left|s_X(e) - t^{k-2} f'(X)^{\binom{k}{2}-1} \right| < \mu f'(X)^{\binom{k}{2}-1} t^{k-2}.
$$ \qed
\end{lemma}
\begin{lemma}\label{l:32}
With probability at least $1-1/n$,
$$
|E(Z_X)|  > (1-2\zeta)\binom{k}{2}\frac{n^2}{m^2}f'(X).
$$ \qed
\end{lemma}

Since $R$ contains at most $O(m^k)$ copies of $K_k$,
we have that with probability at least $1-O(m^k/n) - O(m^{k+3}/n) > 0$ (here we assume again that $n_0$ is sufficiently large) {\em all}
copies $X$ of $K_k$ in $R$ with $f'(X) > m^{1-k}$
satisfy the statements of Lemma \ref{l:31} and Lemma \ref{l:32}.
We therefore fix a partition $\Q$ for which Lemma \ref{l:31} and Lemma \ref{l:32} hold
for all such $X$.

Let $H \in \C_k$.
Let $X$ be a copy of $K_k$ in $R$ with $f'(X) > m^{1-k}$ that is partite-color-isomorphic to $H$.
We construct an $r$-uniform hypergraph $L_X$ as follows.
The vertices of $L_X$ are the edges of $Z_X$.
The edges of $L_X$ correspond to the edge sets of the subgraphs of $Z_X$ that are partite-color-isomorphic
to $H$. We claim that this hypergraph satisfies the conditions of Lemma \ref{l:23}.
Indeed, let $q$ denote the number of vertices of $L_X$.
Let $d=t^{k-2} f'(X)^{\binom{k}{2}-1}$. Notice that by Lemma \ref{l:31} {\em all} vertices of $L_X$ have their degrees between
$(1-\mu)d$ and $(1+\mu)d$. Also notice that the co-degree of any two vertices of $L_X$ is at most $t^{k-3}$
as two edges cannot belong, together, to more than $t^{k-3}$ subgraphs of $L_X$ that are partite-color-isomorphic to
$H$. Also observe that for $n_0$ sufficiently large,
$\mu d > t^{k-3}$. By Lemma \ref{l:23} we have a set $\S_X$ of at least $(q/\binom{k}{2})(1-\beta)$ pairwise edge-disjoint subgraphs of $Z_X$
that are partite-color-isomorphic to $H$.
In particular, by Lemma \ref{l:32},
$$
|\S_X| \ge (1-\beta)(1-2\zeta)\frac{n^2}{m^2}f'(X) > (1-2\beta)f'(X)\frac{n^2}{m^2}\;.
$$
Now, let $\X_H$ be the set of all subgraphs of $R$ that are partite-color-isomorphic to $H$.
By definition, $f'(R,H)=\sum_{X \in \X_H}f'(X)$.
Sine trivially $|\X_H| \le m^k$, the total contribution of the elements $X \in \X_H$ with $f'(X) \le m^{1-k}$ to the
sum is at most $m$. Hence,
\begin{eqnarray*}
\left| \cup_{X \in \X_H\;,\;f'(X) > m^{1-k}} \S_X \right| & \ge & (1-2\beta)\frac{n^2}{m^2}\sum_{X \in \X_H\;,\;f'(X) > m^{1-k}}f'(X)\\
& \ge & (1-2\beta)\frac{n^2}{m^2} \left(f'(R,H)-m\right)\\
& = & (1-2\beta)\frac{n^2}{m^2} \left(f^{*}(G^*,H)\frac{m^2}{n^2}-m\right)\\
& = & (1-2\beta)f^{*}(G^*,H)-(1-2\beta)\frac{n^2}{m}\\
& \ge & (1-2\beta)\left(f(G,H) -\frac{\epsilon}{10} n^2 \right)-(1-2\beta)\frac{n^2}{m}\\
& \ge & f(G,H) - \epsilon n^2\;.
\end{eqnarray*}
As the $\S_X$ are pairwise disjoint for distinct $X$, we have
obtained a set $P_H$ of induced subgraphs of $G$ that are color-isomorphic to $H$, such that $|P_H| \ge f(G,H) - \epsilon n^2$.
Notice further that for distinct $H \in \C_k$, the corresponding sets $P_H$ are disjoint.
\qed

\subsection{Tournaments}

The proof of the tournament case of Theorem \ref{t:main-1} is almost identical to the proof of the edge-colored case presented in this section.
One just needs to prove the following analogue of Lemma \ref{l:1} which is the following Lemma \ref{l:1-tour}.
Recall that an {\em orientation} is a directed simple graph without cycles of length $2$.
A $K_k$-subgraph of an orientation $G$ is a $k$-vertex tournament subgraph of $G$. We similarly define a $K_k$-decomposition and a fractional $K_k$-decomposition
of an orientation.
\begin{lemma}\label{l:1-tour}
Let $k \ge 3$ be an integer and let $\epsilon > 0$. There exists $n_0=n_0(k,\epsilon)$ such that the following holds.
Suppose $G$ is an orientation with $n > n_0$ vertices.
Let $f$ be a fractional $K_k$-decomposition of $G$. Then for every $H \in \T_k$ there is a set $P_H$ of induced subgraphs of $G$ that are
isomorphic to $H$, such that $|P_H| \ge f(G,H) - \epsilon n^2$. Furthermore any two elements of $P=\cup_{H \in \T_k}P_H$ intersect in at most one vertex.
\end{lemma}
The proof of Lemma \ref{l:1-tour} is completely analogous to the proof of Lemma \ref{l:1-tour} where instead of using the 
colored version of Szemer\'edi's regularity lemma (Lemma \ref{l:21}) we use the directed version of the lemma.
We refer to \cite{NY-2007} which contains this directed version of the main result of \cite{yuster-2005} and therefore implies Lemma \ref{l:1-tour}.
We therefore obtain the following corollary, whose proof is analogous to that of corollary \ref{c:1}.

\begin{corollary}\label{c:1-tour}
Let $k \ge 3$ be an integer, let ${\bf v} \in {\mathbb R}^{|\T_k|}$, and let $\gamma > 0$. 
There exists $N_{\ref{c:1-tour}}=N_{\ref{c:1-tour}}(k,\gamma,{\bf v})$ such that the following holds for all orientations $G$ with $n > N_{\ref{c:1-tour}}$ vertices
which have a fractional $K_k$-decomposition.
For every $H \in \T_k$ there is a set $P_H$ of induced subgraphs of $G$ that are isomorphic to $H$, such that any two elements of $P=\cup_{H \in \T_k}P_H$ intersect in  at most one vertex. Furthermore,
\begin{eqnarray*}
(a) & |P| \ge \frac{|E(G)|-\gamma n^2}{\binom{k}{2}}\;.\\
(b) & \sum_{H \in \T_k} {\bf v_H}|P_H| \ge \frac{|E(G)|}{\binom{k}{2}}\nu^*_{\bf v}(G)-\gamma n^2\;.
\end{eqnarray*}
\end{corollary}
Finally, we need the analogue of Lemma \ref{l:main} for the tournament setting. The lemma is proved in exactly the same way using Lemmas \ref{l:turan} and \ref{l:2}
(which stay intact; recall that they do not depend on the setting, whether it is tournaments or edge colored graphs) and using the straightforward Lemma \ref{l:3}
(whose statement stays intact, but in its proof $G$ is a tournament instead of a $C$-colored complete graph). We therefore obtain the following lemma which immediately implies the second part of Theorem \ref{t:main-1}.
\begin{lemma}\label{l:main-tour}
Let $k \ge 3$ be an integer, let ${\bf v} \in {\mathbb R}^{|\T_k|}$ be indexed by $\T_k$, and let $\epsilon > 0$.
Then there exist $N_{\ref{l:main-tour}}=N_{\ref{l:main-tour}}(\epsilon,{\bf v})$ such that the following holds.
Let $G$ be a tournament which is $k$-divisible and with $n > N_{\ref{l:main-tour}}$ vertices.
Then, $G$ is $K_k$-decomposable and, furthermore $\nu_{\bf v}(G) \ge (\nu^*_{\bf v}-\epsilon)(1 - \epsilon)-\epsilon$.
\end{lemma}

\section{Applications}

\subsection{Triangles in tournaments}

Our first application of Theorem \ref{t:main-1} concerns the simplest case $k=3$ for tournaments. Note that here we have $\T_3=\{T_3,C_3\}$ where $T_3$ denotes the transitive triangle and $C_3$ denoted the directed (cyclic) triangle. Recall from the introduction that the various possibilities for $\nu_{\bf v}$ for ${\bf v} \in {\mathbb R}^{|{\T_3}|}$ reduce to the cases where the smallest coordinate of ${\bf v}$ is $0$ and the largest coordinate is $1$ (if ${\bf v}$ is the constant vector, then trivially
$\nu_{\bf v}$ equals that constant). We furthermore see that the case of ${\bf v}(T_3)=0$ is trivial since the sequence of transitive $n$-vertex tournaments
shows that $\nu_{\bf v}=0$ in this case. Hence the {\em only} vector for which $\nu_{\bf v}$ is nontrivial to evaluate is the one which assigns ${\bf v}(T_3)=1$ and ${\bf v}(C_3)=0$.
\begin{conj}
Let ${\bf v} \in {\mathbb R}^{|{\T_3}|}$ where ${\bf v}(T_3)=1$ and ${\bf v}(C_3)=0$. Then $\nu_{\bf v}=1$.
\end{conj}
In \cite{yuster-2004} it was conjectured that every tournament can be packed with $\lceil n(n-1)/6 - n/3 \rceil=(1-o(1))\binom{n}{2}/\binom{3}{2}$ edge-disjoint copies of $T_3$ (and, if true, this conjectured value is shown there to be optimal).
However, notice that even if the conjecture in \cite{yuster-2004} is true, this by no means implies that $\nu_{\bf v}=1$, since we have no guarantee that a very large
$T_3$-packing is part of a triangle decomposition (notice also that a triangle decomposition exists whenever $n \equiv 1,3 \bmod 6$, by Kirkman's Theorem).

Here we prove the following theorem which implies Theorem \ref{t:applic-1}.
\begin{theorem}\label{t:applic-1-lower}
Let ${\bf v} \in {\mathbb R}^{|{\T_3}|}$ where ${\bf v}(T_3)=1$ and ${\bf v}(C_3)=0$. Then, $\nu_{\bf v} \ge \frac{85}{98}$.
In particular, for all $n \equiv 1,3 \bmod 6$, every tournament on $n$ vertices has a triangle decomposition where the number of $C_3$ in the decomposition
is at most $\frac{13}{98 \cdot 6}n^2(1+o_n(1))$.
\end{theorem}
\Proof
By Theorem \ref{t:main-1}, it suffices to prove that $\nu^*_{\bf v} \ge \frac{85}{98}$.
A computer assisted proof (outlined below) yields that $\nu^*_{\bf v}(14)=\frac{78}{91}$.
As proved in Corollary 2.7 in \cite{KY-2008} \footnote{That corollary is used in \cite{KY-2008} for fractional triangle packings but it is identical for fractional triangle decompositions.} following an iterative improvement argument appearing first in Lemma 2.2 of \cite{KS-2004}, 
$\nu^*_{\bf v}(n) \ge (\nu^*_{\bf v}(r)(r-1)+1)/r - o_n(1)$. So, plugging in the case $r=14$, $\nu^*_{\bf v}(14)=\frac{78}{91}$
and taking the limit yields $\nu^*_{\bf v} \ge \frac{85}{98}$.

So, it remains to show that $\nu^*_{\bf v}(14) \ge \frac{78}{91}$.
Let us recall that this means that for {\em every} tournament $G$ on $14$ vertices, there is a fractional triangle decomposition $f$ such that
$\nu^*_{\bf v}(f) \ge \frac{78}{91}$, or, equivalently, that $D_{\bf v}(f) \ge \frac{78}{91} \cdot \binom{14}{2}/\binom{3}{2} = 26$.
In our case, since ${\bf v}(T_3)=1$ and ${\bf v}(C_3)=0$, this means that the sum of the values of $f$ on all $T_3$ copies of $G$
is at least $26$. As a side note, we observe that any $14$ vertex tournament $G_0$ that is obtained by taking three disjoint sets of vertices $A,B,C$
with $|A|=|B|=5$ and $|C|=4$ and orienting all edges from $A$ to $B$, from $B$ to $C$ and from $C$ to $A$ (the orientations of edges with both endpoints in the same part
is arbitrary), has the property that each of its $T_3$
copies contains an edge with both endpoints in the same part. So, for any fractional triangle decomposition $f$,
the sum of the values of $f$ on all $T_3$ copies of such a $G_0$
is at most the number of edges with both endpoints in the same part which is $\binom{5}{2}+\binom{5}{2}+\binom{4}{2}=26$.
Thus, we always have $\nu^*_{\bf v}(14) \le \frac{78}{91}$.

The naive computational approach would therefore be as follows. Generate all $14$-vertex tournaments $G$ (say, up to isomorphism).
For each such $G$, write down the linear programming problem which finds a fractional triangle decomposition which maximizes the
sum of the values it assigns to the $T_3$ elements of $G$, and verify that this maximum, denoted by $D^*(G)$ is always at least $26$.
This naive approach is infeasible since the number of (pairwise non-isomorphic) tournaments on $14$ vertices is more than any computer can handle
(already the number of $14$-vertex strongly connected tournaments on $14$ vertices is $28304491788158056$ by the OEIS), moreover running a (rather large) linear programming instance on each. Instead we take the following significantly better approach.

We call a tournament $G$ on $r+1$ vertices an {\em extension} of a tournament $G'$ on $r$ vertices, if $G'$ is a subgraph of $G$.
Notice that a tournament on $r$ vertices has at most $2^r$  extensions as can be seen by adding a new vertex and considering all possible
orientations of its $r$ incident edges.
The following simple lemma is immediate from the proof of Lemma \ref{l:3} (the construction of $f$ there).
\begin{lemma}\label{l:extension}
Let $G$ be a tournament with $r+1$ vertices. If $D^*(G) < t$,  then it is an extension of some $G'$ with $D^*(G') < t \cdot \frac{r-1}{r+1}$. \qed
\end{lemma}
For $r \ge 3$ and a real $t$, let $\T_r(t)$ denote the set of all $r$-vertex tournaments $G$ with $D^*(G) < t$.
So, our goal is to prove that $\T_{14}(26) = \emptyset$.
By lemma \ref{l:extension}, it suffices to check all extensions of $\T_{13}(26 \cdot \frac{12}{14}) \subseteq \T_{13}(22.3)$.
In turn, it suffices to check all extensions of $\T_{12}(18.86)$. 
In turn, it suffices to check all extensions of $\T_{11}(15.72)$. 
In turn, it suffices to check all extensions of $\T_{10}(12.86)$. 
So, we start by generating all non-isomorphic tournaments on $10$ vertices. There are known lists of such tournaments, see
\url{https://users.cecs.anu.edu.au/~bdm/data/digraphs.html}.
There are only $9733056$ such tournaments. We denote this set by $M_{10}$.
For each $G \in M_{10}$, we run the corresponding linear program to compute $D^*(G)$. If $D^*(G) \ge 12.86$ then, as shown earlier, we are not worried,
as we are not missing anything by not checking extensions of such $G$.
However, if $D^*(G) < 12.86$, we say that $G$ is below the threshold, so we generate all $2^{10}$ extensions of $G$
(we don't mind generating isomorphic tournaments, as the time required to check isomorphisms would be larger).
Doing it for all $G$ on $10$ vertices which are below the threshold, yields a multiset of tournaments on $11$ vertices, call it $M_{11}$.
Notice that by the above, we know that if some $G$ on $14$ vertices has $D^*(G) < 26$, then it contains an element of $M_{11}$ as a subgraph.
Now, for each tournament $G \in M_{11}$, we run the corresponding linear program. If $D^*(G) < 15.72$, we generate all $2^{11}$ extensions of $G$.
This yields a multiset of tournaments on $12$ vertices, call it $M_{12}$. 
For each tournament $G \in M_{12}$, if $D^*(G) < 18.86$, we generate all $2^{12}$ extensions of $G$.
This yields a multiset of tournaments on $13$ vertices, $M_{13}$.
For each tournament $G \in M_{13}$, if $D^*(G) < 22.3$, we generate all $2^{13}$ extensions of $G$.
This yields a multiset of tournaments on $14$ vertices, $M_{14}$.
Finally, we check all tournaments in $M_{14}$ to verify that $D^*(G) \ge 26$ for each of them.
This procedure is summarized in Table \ref{t:thresholds}. The table also lists for each $r=10,\ldots,14$ the size of the (multi)set $M_r$,
the number of elements of $M_r$ that are below the threshold, which means that $|M_{r+1}|$ is precisely $2^r$ times larger than this amount.
We also list the lowest value of $D^*(G)$ encountered during the search.
\begin{table}
	\centering
		\begin{tabular}{|c||c|c|c|c|}
		 \hline
		  $r$  		& Size of $M_r$ & Threshold value & Below threshold	& Lowest value	\\ \hline
		  $10$ 	& $9733056$ 		& $12.86$ 				& $16$						& $12$					\\ \hline
		  $11$ 	& $16384$ 			& $15.72$ 				& $256$						& $15$					\\ \hline
		  $12$ 	& $524288$ 			& $18.86$ 				& $2048$					&	$18$					\\ \hline
		  $13$ 	& $8388608$ 		& $22.3$ 					& $98304$					& $22$					\\ \hline
		  $14$ 	& $805306368$		& $26$ 						& $0$ 						& $26$					\\ \hline
		\end{tabular}
	\caption{The procedure for verifying that $D^*(G) \ge 26$ for all $14$-vertex tournaments.}
	\label{t:thresholds}
\end{table}

The code of the program that performs the procedure above can be found in
\url{https://github.com/raphaelyuster/Vector-valued-decompositions}.
The program runs in fewer than five days on standard personal computer equipment.
The program uses the well-established linear programming package {\ttfamily lp-solve} which has a very efficient and simple to use api, see
the {\ttfamily lp-solve} package homepage can be found at
\url{https://sourceforge.net/projects/lpsolve}.
Let us note that the linear programming instances are easy to generate. Suppose $G$ is a tournament on $r$ vertices.
Generate a variable for each triple of vertices of $G$, so there are $\binom{r}{3}$ variables.
To compute $D^*(G)$ one should maximize the sum of the variables that correspond to triples that induce a $T_3$.
The constraints are:
For each pair of vertices $i,j$ of $G$, the sum of the variables that correspond to triples that contain both $i,j$ should be precisely $1$.
Hence there are $\binom{r}{2}$ such constraints.
Furthermore, we require that all variables are nonnegative. So there are $\binom{r}{3}$ such constraints.
This completes the proof of Theorem \ref{t:applic-1-lower}.

\subsection{Essentially avoidable graphs}

In order to prove Theorem \ref{t:applic-2} we need to extend the notion of $\nu^*_{\bf v}(f)$ to fractional packings. Formally, for ${\bf v} \in {\mathbb R}^{|\C_k|}$ indexed by $\C_k$
and a fractional $K_k$-packing $f$ of a $C$-colored graph $G$,
let $D_{\bf v}(f) = \sum_{H \in \C_k} {\bf v}_H f(G,H)$. After normalizing we define
$\nu^*_{\bf v}(f) = D_{\bf v}(f)\binom{k}{2}/|E(G)|$.
Since the result and proof in \cite{yuster-2005} applies to fractional packings,
so do Lemma \ref{l:1} and Corollary \ref{c:1}. Restated for fractional packings, Corollary \ref{c:1} becomes:
\begin{corollary}\label{c:2}
Let $C$ be a  finite set of colors, let $k \ge 3$ be an integer, let ${\bf v} \in {\mathbb R}^{|\C_k|}$, and let $\gamma > 0$. 
There exists $N_{\ref{c:2}}=N_{\ref{c:2}}(k,C,\gamma,{\bf v})$ such that the following holds for all $C$-colored graphs $G$ with $n > N_{\ref{c:2}}$ vertices.
Suppose $f$ is a fractional $K_k$-packing of $G$. Then
for every $H \in \C_k$ there is a set $P_H$ of induced subgraphs of $G$ that are color-isomorphic to $H$, such that any two elements of $P=\cup_{H \in \C_k}P_H$ intersect in  at most one vertex. Furthermore,
\begin{eqnarray*}
(a) & |P| \ge \frac{\left(\sum_{H \in \C_k}f(G,H)\right)-\gamma n^2}{\binom{k}{2}}\;.\\
(b) & \sum_{H \in \C_k} {\bf v_H}|P_H| \ge \frac{|E(G)|}{\binom{k}{2}}\nu^*_{\bf v}(f)-\gamma n^2\;.
\end{eqnarray*}
\end{corollary}

We say that a fractional $K_k$-packing $f$ of a $C$-colored graph $G$ is {\em $\delta$-close to a fractional decomposition} if $\sum_{H \in \C_k}f(G,H) \ge |E(G)| - \delta n^2$. We say that a binary vector ${\bf v} \in {\mathbb R}^{|\C_k|}$ is {\em nice} if for every $\delta > 0$, if $n$ is sufficiently large, then
every $C$-colored complete graph $G$ on $n$ vertices has a fractional packing $f$ that is $\delta$-close to fractional decomposition and with $\nu^*_{\bf v}(f) \ge 1-\delta$.
With these definitions, together with Corollary \ref{c:2}, Lemma \ref{l:turan} and Lemma \ref{l:2}, the following lemma is proved in the same way Lemma
\ref{l:main} is proved.
\begin{lemma}\label{l:almost}
Let $C$ be a  finite set of colors, let $k \ge 3$ be an integer, let ${\bf v} \in {\mathbb R}^{|\C_k|}$ be a nice vector indexed by $\C_k$, and let $\epsilon > 0$.
Then there exist $N_{\ref{l:almost}}=N_{\ref{l:almost}}(\epsilon,{\bf v})$ such that the following holds.
Let $G$ be a $C$-colored complete graph which is $k$-divisible and with $n > N_{\ref{l:almost}}$ vertices.
Then, $G$ is $K_k$-decomposable and $\nu_{\bf v}(G) \ge 1-\epsilon$. \qed
\end{lemma}

Let $C=\{red,blue\}$. A binary vector ${\bf v} \in {\mathbb R}^{|\C_k|}$ therefore corresponds to a characteristic vector of a subset $\F^* \subset \F_k$,
where ${\bf v}_H = 1$ if and only if the blue edges of $H$ correspond to a graph form $\F^*$.

\noindent
{\bf Proof of Theorem \ref{t:applic-2}.}
Let $k \ge 5$ be odd and $\F^* \subset \F_k$ be the family of all graphs $H$ on $k$ vertices such that both $H$ and its complement are
Eulerian. Let ${\bf v}$ be the corresponding characteristic vector of $\F_k \setminus \F^*$.
Theorem 3 of \cite{yuster-2019} implies that ${\bf v}$ is nice.
By Lemma \ref{l:almost}, if $G$ is a $k$-divisible graph on $n > N_{\ref{l:almost}}$ vertices, then
$\nu_v(G) \ge 1-\epsilon$. Thus, $\F^*$ is essentially avoidable.

For a graph $H$ on $k$ vertices, let ${\bf v}$ be the characteristic vector of $\F_k \setminus \{H\}$.
Let ${\cal U}_k \subseteq \F_k$ be the set of graphs $H$ whose corresponding characteristic vector of $\F_k \setminus \{H\}$ is not nice.
By Lemma \ref{l:almost}, this implies that each $H \notin {\cal U}_k$ is essentially avoidable. 
Theorem 2 of \cite{yuster-2019} implies that $|{\cal U}_k|=o(|\F_k|)$. Hence the second part of the theorem follows. \qed

\section{Small \texorpdfstring{$k$}{k}}

We start this section by considering $\nu_{\bf v}$ for ${\bf v} \in {\mathbb R}^{|\F_3|}$ which is the first nontrivial case.
By Theorem \ref{t:main-1} it suffices to determine $\nu^*_{\bf v}$ and as noted in the introduction our problem is reduced to vectors whose smallest
coordinate is $0$ and whose largest coordinate is $1$. We call such vectors {\em normalized}.

Observe that $\F_3=\{K_3,P_3,Q_3,I_3\}$ where $P_3$ denotes the path on three vertices, $I_3=K_3^c$ is the independent set on $3$ vertices
and $Q_3=P_3^c$. We will use the convention of writing ${\bf v} = ({\bf v}(K_3),{\bf v}(P_3),{\bf v}(Q_3),{\bf v}(I_3))$.
The following proposition determines $\nu_{\bf v}$ for a significant amount of normalized vectors, which include in particular all binary vectors.
\begin{prop}\label{p:1}
Let ${\bf v} \in {\mathbb R}^{|\F_3|}$ be a normalized vector.
\begin{enumerate}
\item
$\frac{1}{4}\min \{{\bf v}(K_3)\, ,\, {\bf v}(I_3)\} \le  \nu_{\bf v} \le \min \{{\bf v}(K_3)\, ,\, {\bf v}(I_3)\}$.
\item
If ${\bf v}=(1,0,\beta,\alpha)$ or ${\bf v}=(\alpha,\beta,0,1)$,  then $\nu_{\bf v} = \frac{\alpha}{4}$.
\end{enumerate}
\end{prop}
\Proof
For the first part of the proposition, consider first $G=K_n$. Here each $3$-vertex subgraph is a $K_3$ so we obtain $\nu^*_{\bf v}(G)={\bf v}(K_3)$.
Similarly, for $G=I_n$ we have $\nu^*_{\bf v}(G)={\bf v}(I_3)$. Hence, $\nu^*_{\bf v}(n) \le \min \{{\bf v}(K_3)\, ,\, {\bf v}(I_3)\}$
so $\nu_{\bf v} = \nu^*_{\bf v} \le \min \{{\bf v}(K_3)\, ,\, {\bf v}(I_3)\}$.

We recall a theorem of Goodman \cite{goodman-1959} who proved that in any $n$-vertex graph,
$\frac{1}{4}\binom{n}{3}(1-o_n(1))$ of the sets of $3$ vertices induce either a $K_3$ or an $I_3$. Hence, the fractional decomposition $f$ which
assigns a value of $1/(n-2)$ to each $3$-set of vertices has $\nu^*_{\bf v}(f) \ge (1-o_n(1))\frac{1}{4}\{ \min \{{\bf v}(K_3)\, ,\, {\bf v}(I_3)\}$.
Hence, $\nu_{\bf v} = \nu^*_{\bf v} \ge \frac{1}{4}\{ \min \{{\bf v}(K_3)\, ,\, {\bf v}(I_3)\}$.

For the second part of the proposition, note first that for ${\bf v}=(1,0,\beta,\alpha)$ or ${\bf v}=(\alpha,\beta,0,1)$,
the aforementioned lower bound implies that $\nu_{\bf v} = \nu^*_{\bf v} \ge \frac{\alpha}{4}$.

For the upper bound, we consider first the case ${\bf v}=(1,0,\beta,\alpha)$.
Let $G$ be the complete balanced bipartite graph on $n$ vertices, where the sides are $A,B$ with $|A|=\lfloor n/2 \rfloor$ and $|B|=\lceil n/2 \rceil$.
We assume that $n \equiv 1,3 \bmod 6$ so that there is a $3$-decomposition of $G$. Notice that this implies that $n$ is odd and that $|A||B|$ is even.
Consider some $3$-decomposition $L$ of $G$. As the $|A||B|=(n^2-1)/4$ edges of $G$ must be packed, there exist precisely $|A||B|/2 = (n^2-1)/8$ elements of $L$
that are isomorphic to $P_3$. These elements also contain $(n^2-1)/8$ pairs  of vertices with both endpoints in the same part, so $L$ has precisely
$n(n-1)/6-(n^2-1)/8=(n^2-4n+3)/24$ elements isomorphic to $I_3$. This proves that
$$
\nu_{\bf v}(n) \le \alpha \frac{n^2-4n+3}{24} \cdot \frac{3}{\binom{n}{2}} = \frac{\alpha}{4}\cdot \frac{n-3}{n}
$$
proving that $\nu_{\bf v} \le \frac{\alpha}{4}$. The case ${\bf v}=(\alpha,\beta,0,1)$ is proved analogously by taking complements.
\qed

Observe that Proposition \ref{p:1} determines $\nu_{\bf v}$ for all binary vectors. It is always $0$ unless 
${\bf v}(K_3)={\bf v}(I_3)=1$ in which case it is $\frac{1}{4}$ except for the trivial case ${\bf v}=(1,1,1,1)$ where we have $\nu_{\bf v}=1$.
Still, Proposition \ref{p:1} does not cover all possible normalized vectors, so we raise the following problem.
\begin{prob}
Determine $\nu_{\bf v}$ for all (normalized vectors) ${\bf v} \in {\mathbb R}^{|\F_3|}$.
\end{prob}

Moving to the next case $k=4$, we do not know the value of $\nu_{\bf v}$ even for all binary vectors.
Notice that $|\F_4|=11$ as there are $11$ distinct $4$-vertex graphs. While the all-1 vector trivially has $\nu_{\bf v}=1$, we do not know of a single
normalized vector for which $\nu_{\bf v} = 1$. So, a realistic open problem is the following.
\begin{prob}
Determine the normalized vectors ${\bf v} \in {\mathbb R}^{|\F_4|}$ for which $\nu_{\bf v}=1$.
\end{prob}
Small examples suggest that it is plausible that the binary vector which assigns $1$ to all graphs in $\F_4$ except $C_4$ and assigns $0$ to $C_4$ has $\nu_{\bf v}=1$.
Finally, note that for $k=5$ we know of a normalized binary vector which has $\nu_{\bf v}=1$. Indeed by Theorem \ref{t:applic-2}, the vector which assigns
$1$ to all graphs in $\F_5$ except $C_5$ and assigns $0$ to $C_5$ has $\nu_{\bf v}=1$.

\section{The random graph}

In this section we asymptotically determine $\nu_{\bf v}(G)$ for almost all $k$-decomposable graphs $G$ and for all ${\bf v} \in {\mathbb R}^{|\F_k|}$.
The main result in this section is stated for the Erd\H{o}s-R\'enyi random graph probability space $\G(n,p)$ where $0 < p < 1$ is a constant.
Recall that a property which holds for almost all $G \sim \G(n,\frac{1}{2})$ is referred to as a property that holds for almost all graphs.

Recall that an $n$-vertex graph $G \sim \G(n,p)$
is obtained by independently deciding for each pair of vertices whether it is an edge with probability $p$.
Now, suppose $n$ is such that graphs with $n$ vertices are $k$-decomposable (recall that by Wilson's Theorem this holds for all $n$ sufficiently large
such that $K_n$ is $K_k$-divisible). Then, for $G \sim \G(n,p)$ we have that $\nu_{\bf v}(G)$ is a random variable, and hence our ultimate goal would
be to show that $\nu_{\bf v}(G)$ converges in distribution to a constant, and determine this constant.
Indeed this is the main result in this section.

To define the constant to which $\nu_{\bf v}(G)$ converges in distribution, we set up a small (constant size) linear program.
Let ${\bf v} \in {\mathbb R}^{|\F_k|}$ and consider the linear program $LP({\bf v},p)$ defined as follows.
\begin{eqnarray*}
& {\bf max} & \sum_{H \in \F_k} {\bf v}_H x_H \\
& {\bf s.t.} & \sum_{H \in \F_k} \left(e(H)-p\binom{k}{2}\right)x_H = 0\, ,\\
& & \sum_{H \in \F_k} x_H = 1\, ,\\
& & x_H \ge 0 ~~~ \forall H \in \F_k\, .
\end{eqnarray*}
Clearly $LP({\bf v},p)$ is feasible since setting $x_{K_k}=p$ and $x_{I_k}=(1-p)$ and setting all other variables to $0$, all constraints are satisfied.
Therefore, let $s({\bf v},p)$ denote the optimal solution of $LP({\bf v},p)$.
Our main theorem follows. Notice that when we write $n \rightarrow \infty$ we only consider $n$ such that graphs with $n$ vertices are $k$-decomposable.
\begin{theorem}\label{t:random}
Let ${\bf v} \in {\mathbb R}^{|\F_k|}$ and let $0 < p < 1$. For every $\epsilon > 0$, $G \sim \G(n,p)$ satisfies
$$
\lim_{n \rightarrow \infty}  \Pr \left[ |\nu_{\bf v}(G) - s({\bf v},p)| < \epsilon \right] = 1 \;.
$$
\end{theorem}
\Proof
We first prove that $\Pr \left[ \nu_{\bf v}(G) \ge s({\bf v},p)+ \epsilon \right] = o_n(1)$.
To this end, we don't even need to assume that $G$ is a random graph. All that suffices is to assume that $e(G) = p\binom{n}{2} \pm o(n^2)$,
which trivially holds with probability $1-o_n(1)$ for $G \sim \G(n,p)$.
So, assume that $G$ is an $n$-vertex, $k$-decomposable graph with $e(G) = p\binom{n}{2} \pm o(n^2)$. We prove that $\nu_{\bf v}(G) \le s({\bf v},p)+ \epsilon$.
Take an optimal $k$-decomposition $L$ of $G$ with respect to ${\bf v}$.
For $H \in \F_k$, let $L_H$ be the subset of $L$ whose elements are isomorphic to $H$
and let $y_H = |L_H|/|L|$. Observe that $y_H \ge 0$ and that $\sum_{H \in \F_k} y_H=1$.
Next, observe that $e(H)|L_H|$ is the total number of edges of $G$ in all the elements of  $L_H$, and since $L$ is a decomposition,
we have that $\sum_{H \in \F_k} e(H)|L_H|=e(G)$ and so $\sum_{H \in \F_k} e(H)y_H=e(G)/|L|$.
But since $|L|=\binom{n}{2}/\binom{k}{2}$ and since $e(G) = p\binom{n}{2} \pm o(n^2)$ we have that
$\sum_{H \in \F_k} \left(e(H)-p\binom{k}{2}\right)y_H = o_n(1)$.
Hence, there exist $z_H$ such that $|z_H-y_H| = o_n(1)$ for all $H \in \F_k$ such that $z_H$ form a feasible solution of $LP({\bf v},p)$
and such that for all $n$ sufficiently large, $\sum_{H \in \F_k} {\bf v}_H(y_H-z_H) \le \epsilon$.
As the $z_H$ form a feasible solution we get that
$$
\nu_{\bf v}(G) = \nu_{\bf v}(L) = \sum_{H \in \F_k} {\bf v}_H y_H  \le \epsilon + \sum_{H \in \F_k} {\bf v}_H z_H \le s({\bf v},p)+ \epsilon\;.
$$

We next prove that $\Pr \left[ \nu_{\bf v}(G) \ge s({\bf v},p) - \epsilon \right] = 1-o_n(1)$.
We will assume for simplicity that $p$ is rational. This can be assumed since for given $n,\epsilon, {\bf v}$, the function
$\Pr \left[ \nu_{\bf v}(G) \ge s({\bf v},p) - \epsilon \right]$ where $G \sim \G(n,p)$ is continuous in $p$.

Now, as $p$ is rational, so is $s({\bf v},p)$ and there is an optimal solution ${\bf x}= \{x_H:H \in \F_k\}$ where all the $x_H$ are rational.
By taking a common denominator $d$, we denote $x_H=a_H/d$ where the $a_H$ are nonnegative integers not exceeding $d$.
It will be convenient to view $G \sim \G(n,p)$ as an edge colored $K_n$ where the blue edges are the edges of $G$ and the red edges are
the non-edges of $G$ and similarly view the elements of $\F_K$ as blue-red edge colored $K_k$.

We construct a gadget blue-red edge-colored graph $D$ as follows. $D$ consists of $d$ edge-disjoint copies of $K_k$ (any such graph $D$ suffices).
For each $H \in F_k$ precisely $a_H$ of the $K_k$ comprising $D$ are color-isomorphic to $H$. Notice that $D$ has precisely $d\binom{k}{2}$ edges,
where $\sum_{H \in \F_k} a_He(H)$ of them are colored blue and the others are colored red.
But observe that since the $x_H$ form a feasible solution to $LP({\bf v},p)$, this also means that the number of blue edges of $D$ is
$pd\binom{k}{2}$ and the number of red edges is $(1-p)d\binom{k}{2}$.

Let $r > 1$ be the smallest integer such that $K_r$ has a $D$-decomposition. By Wilson's Theorem, $r$ exists.
Let $R$ be a blue-red edge coloring of $K_r$ obtained by taking a $D$-decomposition of $K_r$, and coloring each element of this decomposition
such that it is color isomorphic to $D$. Observe that the number of blue edges of $R$ is $p\binom{r}{2}$ and the number of red edges is
$(1-p)\binom{r}{2}$.

Now we consider $G \sim \G(n,p)$ (recall that $G$ is viewed as a blue-red edge-colored $K_n$). We construct an $\binom{r}{2}$ uniform hypergraph $M$ as follows.
The vertices of $M$ are the $\binom{n}{2}$ edges of $G$. The edges of $M$ are all the $K_r$-subgraphs of $G$ that are color-isomorphic to $R$.
We observe some properties of $M$ which stem from the fact that $G \sim \G(n,p)$.
What is the degree of a blue vertex of $M$, or, stated equivalently, what is the number of copies of $R$ in $G$ that contain a given blue edge?
For an $r$-set of vertices of $G$, let $q$ denote the probability that it induces $R$. It doesn't really matter what $q$ is, but nevertheless it is easy
to compute it: $q=p^{p\binom{r}{2}}(1-p)^{(1-p)\binom{r}{2}}r!/aut(R)$ where $aut(R)$ is cardinality of the color-preserving automorphism group of $R$.
For a given pair of vertices $u,v$ and for an additional set $W$ of $r-2$ vertices, what is the probability that $W \cup \{u,v\}$ induces $R$
and that $(u,v)$ is blue? Since only a $p$ fraction of edges of $R$ are blue, and given that $W \cup \{u,v\}$ induces $R$, $(u,v)$ is equally likely
to be any edge of $R$, the probability that $W \cup \{u,v\}$ induces $R$ and that $(u,v)$ is blue is precisely $pq$.
Now, given that $(u,v)$ is blue, the probability of an additional subset $W$ of $r-2$ vertices to induce together with $u,v$ a copy of $R$
is, by conditional expectation, precisely $pq/p=q$. Hence, the expected degree of a blue vertex of $M$ is precisely $q\binom{n-2}{r-2}$.
Similarly given that $(u,v)$ is red, the probability of an additional subset $W$ of $r-2$ vertices to induce together with $u,v$ a copy of $R$
is, by conditional expectation, precisely $(1-p)q/(1-p)=q$ so the the expected degree of a red vertex of $M$ is also precisely $q\binom{n-2}{r-2}$.
Since the degree of a vertex of $M$ (i.e. edge of $G$) is a random variable which is the sum of $\binom{n-2}{r-2}$ indicator random variables
and each variable only depends on $O(n^{r-3})$ other variables, we have by Janson's inequality that for all $n$ sufficiently large,
the probability that all vertices of $M$ have their degrees $q\binom{n-2}{r-2} + o(n^{r-2})$ is $1-o_n(1)$.
Another (trivial) property of $M$ is that the co-degree of any two vertices of $M$, or equivalently, the number of copies of $R$ in $G$ that contain two distinct given edges
is $O(n^{r-3})$.

Given these properties of $M$ we can now apply Lemma \ref{l:23} (the Frankl-R\"odl hypergraph matching theorem) which states that
with probability $1-o_n(1)$, $M$ has a matching covering all but $o(V(M))$ of the vertices of $M$. In other words, with probability $1-o(1)$, there is
a packing of $G$ with pairwise edge-disjoint copies of $R$, such that the number of unpacked edges is $o(n^2)$.
But now recall that each copy of $R$ decomposes into $D$ and each copy of $D$ contains, for each $H \in F_k$, precisely $a_H$ pairwise edge-disjoint
$K_k$ subgraph that are color-isomorphic to $H$. But since the $x_H$ are an optimal solution to $LP({\bf v},p)$, we get that
with probability $1-o_n(1)$, there is a $k$-packing $P$ of $G$ such that $\nu_{\bf v}(P) \ge s({\bf v},p) - o_n(1)$.
We can now only slightly modify $P$ to obtain a $k$-decomposition $L$ using Lemma \ref{l:turan} precisely in the same way shown in Lemma \ref{l:main}
where $\nu_{\bf v}(L) \ge \nu_{\bf v}(P) -o_n(1)$. Thus, $\nu_{\bf v}(L) \ge s({\bf v},p) - o_n(1)$ with probability $1-o_n(1)$,
implying that for every $\epsilon > 0$, $\Pr \left[ \nu_{\bf v}(G) \ge s({\bf v},p) - \epsilon \right] = 1-o_n(1)$.

Combining now the two parts of the proof we obtain that $\Pr \left[ |\nu_{\bf v}(G) - s({\bf v},p)| < \epsilon \right] = 1-o_n(1)$, implying the theorem.
\qed

Since $s({\bf v},p)$ can be solved in constant time for every ${\bf v} \in {\mathbb R}^{|\F_k|}$, we can view Theorem \ref{t:random} as saying
that the asymptotic value of $\nu_{\bf v}(G)$ is determined for almost all graphs (using $p=\frac{1}{2}$).

We end this section with an example of a nontrivial case already for $k=3$.
Using the notation of the previous section, consider the vector ${\bf v} \in {\mathbb R}^{|\F_3|}$ defined by ${\bf v}(K_3)=1$, 
${\bf v}(P_3)=\frac{1}{2}$, ${\bf v}(Q_3)=\frac{1}{2}$, ${\bf v}(I_3)=0$ and assume $p=\frac{1}{2}$.
Putting $x_3=x_{K_3}$, $x_2=x_{P_3}$, $x_1=x_{Q_3}$, $x_0=x_{I_3}$, the linear program $LP({\bf v},\frac{1}{2})$ becomes:
\begin{eqnarray*}
& {\bf max} & \frac{1}{2}x_1+\frac{1}{2}x_2+x_3\\
& {\bf s.t.} & -\frac{3}{2}x_0 - \frac{1}{2}x_1 + \frac{1}{2}x_2 + \frac{3}{2}x_3 = 0\, ,\\
& & x_0+x_1+x_2+x_3 = 1\, ,\\
& & x_i \ge 0 ~~~ \forall i \in \{0,1,2,3\}\, .
\end{eqnarray*}
The optimal solution here is $s({\bf v},\frac{1}{2})=\frac{5}{8}$ with $x_1= \frac{3}{4}$, $x_3=\frac{1}{4}$, $x_0=x_2=0$.
Mimicking the proof of Theorem \ref{t:random}, we construct a gadget blue-red edge colored graph $D$ consisting of four edge disjoint triangles.
One triangle is completely blue (this corresponds to one copy of $K_3$), the other three triangles each have two red edges and one blue edge
(this corresponds to three copies of $Q_3$). We observe that $D$ has $12$ edges, $6$ of which are blue and $6$ are red.
As in the proof of Theorem \ref{t:random}, a random graph $G \sim \G(n,\frac{1}{2})$ where the non-edges are colored red and the edges are colored blue
almost surely almost decomposes to $D$. So, as in the theorem, this implies that we have a decomposition $L$ of $G$ into triangles where the number of blue triangles
is roughly $n(n-1)/24$ and the number of triangles with two red edges and one blue edge is roughly $n(n-1)/8$. This yields
that $\nu_{\bf v}(L) = \frac{5}{8}(1-o_n(1))$.

\end{document}